\documentclass[12pt]{amsart}
\usepackage{graphicx}
\title[Skolem-Mahler-Lech Theorem in characteristic $>0$]{A Skolem-Mahler-Lech Theorem in 
Positive Characteristic and Finite Automata}
\author{Harm Derksen}
\thanks{The author was partially supported by NSF grant, DMS 0349019}
\newcommand{\letters}{\renewcommand{\theenumi}{\alph{enumi}}}

\newtheorem{theorem}{Theorem}[section]
\newtheorem{proposition}[theorem]{Proposition}
\newtheorem{lemma}[theorem]{Lemma}
\newtheorem{corollary}[theorem]{Corollary}
\newtheorem{conjecture}[theorem]{Conjecture}
\newcommand{\comp}{\operatorname{comp}}

\theoremstyle{definition}
\newtheorem{definition}[theorem]{Definition}

\newtheorem{example}[theorem]{Example}

\newcommand{\id}{\operatorname{id}}
\newcommand{\Ann}{\operatorname{Ann}}
\newcommand{\Res}{\operatorname{Res}}
\newcommand{\N}{{\mathbb N}}
\newcommand{\Z}{{\mathbb Z}}
\newcommand{\F}{{\mathbb F}}
\newcommand{\Q}{{\mathbb Q}}
\newcommand{\R}{{\mathbb R}}

\begin{document}
\maketitle
\begin{abstract}
Lech proved in 1953 that the set of zeroes of a linear recurrence sequence
in a field of characteristic 0 is the union
of a finite set and finitely many infinite arithmetic progressions.
This result is known as the Skolem-Mahler-Lech theorem. Lech gave a
counterexample to a similar statement in positive characteristic.
We will present some more pathological examples.
We will state and prove a correct analog of the Skolem-Mahler-Lech theorem
in positive characteristic. The zeroes of a recurrence sequence in positive
characteristic can be described using finite automata.
\end{abstract}
\tableofcontents
\section{Introduction}
Suppose that $R$ is a commutative ring (with 1) and $M$ is a (left) $R$-module.
An {\it infinite $M$-sequence} is an element in $M^{\N}$,
where $\N=\{0,1,2,\dots\}$ is the set of nonnegative integers. We say
that $a\in M^{\N}$ satisfies an {\it $R$-recurrence relation of order $d$}
if there exist $\gamma_0,\gamma_1,\dots,\gamma_{d-1}\in R$
such that
\begin{equation}\label{eqrec}
a(n+d)=\gamma_{d-1} a(n+d-1)+\gamma_{d-2}
a(n+d-2)+\cdots+\gamma_0a(n)
\end{equation}
for all $n\in \N$. We will call such a sequence an {\it $R$-recurrence sequence}.
The smallest nonnegative integer $d$ for which $a\in M^{\N}$ satisfies
a recurrence relation of the form (\ref{eqrec}) is called
the {\it order} of the recurrence sequence $a$. For 
 a sequence $a\in M^{\N}$ we define its set of zeroes by
$$
{\mathcal Z}(a)=\{n\in \N\mid a(n)=0\}.
$$

An infinite arithmetic progression is a set of the form $m+n\N$
where $m\in\N$ and $n$ a positive integer. The following result
is the celebrated Skolem-Mahler-Lech theorem.
\begin{theorem}\label{theoSML}
Suppose that $K$ is a field of characteristic $0$ and $a\in K^{\N}$
is a $K$-recurrence sequence. Then ${\mathcal Z}(a)$ is the
union of a finite set and finitely many infinite arithmetic
progressions.
\end{theorem}
This theorem was proved by Skolem (\cite{Skolem}) for $K=\Q$ (the rational numbers),  in 1934, 
by Mahler (\cite{Mahler1}) in 1935
 for $K=\overline{\Q}$ (the algebraic numbers) 
and by Lech for arbitrary fields of characteristic 0 (\cite{Lech,Mahler2}) in
1953.
See also~\cite[\S2.1]{GAIT}. All proofs use an embedding
of $K$ into the $p$-adic completion $\Q_p$ of $\Q$. 

It is possible to bound the number of arithmetic progressions,
and the size of the finite set in
Theorem~\ref{theoSML}
(see~\cite{Schmidt1,Schmidt2} and \cite[Theorem 1.2]{ESS}). Nevertheless 
it is still an open problem
whether ${\mathcal Z}(a)$ can
always be determined for a given $K$-recurrence sequence $a\in K^{\N}$
where $K$ is a field of characteristic 0. In particular, it
is not known if it is decidable whether ${\mathcal Z}(a)=\emptyset$.

The Skolem-Mahler-Lech theorem can be slightly generalized as follows:
\begin{theorem}\label{thm1x}
Suppose that $R$ is a $\Q$-algebra, $M$ is a left $R$-module
and $a\in M^{\N}$ is an $R$-recurrence sequence.
 Then ${\mathcal Z}(a)$ is the
union of a finite set and finitely many infinite arithmetic
progressions.
\end{theorem}
In this paper we will focus on sequences in
fields of positive characteristic. For a prime power $q$ we denote
the field with $q$ elements by $\F_q$. We also define $\F_0=\Q$.
Let $p$ be a prime number. It was  noted
by Lech (\cite{Lech}) that the Skolem-Mahler-Lech theorem without
any modifications is false in positive characteristic:
\begin{example}\label{ex2}
The sequence $a\in \F_p(x)^{\N}$ defined by
$$
a(n)=(x+1)^n-x^n-1
$$
is an $\F_p(x)$-recurrence sequence. The sequence satisfies
$$
a(n+3)=(2x+2)a(n+2)-(x^2+3x+1)a(n+1)+(x^2+x)a(n)
$$
for all $n\in \N$.
The zero set
$$
{\mathcal Z}(a)=\{p^n\mid n\in \N\}
$$
is clearly not the union of a finite set and a finite number
of arithmetic progressions.
\end{example}
Examples such as Example~\ref{ex2} do not yet reveal all the pathologies
that appear in positive characteristic. The following example is new
and stranger.
\begin{example}\label{ex2a}
Define $a\in \F_p(x,y,z)^{\N}$ by
$$
a(n)=(x+y+z)^n-(x+y)^n-(x+z)^n-(y+z)^n+x^n+y^n+z^n.
$$
We have (see Proposition~\ref{prop1})
$$
{\mathcal Z}(a)=\{p^n\mid n\in \N\}\cup \{p^n+p^m\mid n,m\in \N\}.
$$
\end{example}
In \cite{Masser}, Masser gave similar examples of
what he calls ``nested Frobenius type solutions'' to linear equations
over groups in positive characteristic (also called $S$-unit equations).

In order to describe the zero sets of linear recurrence sequences
in positive characteristic, we  make the following definition.
\begin{definition}
Let $p$ be a prime number and $q=p^e$ for some positive integer $e$.
Suppose that $d\geq 1$,
$c_0,c_1,\dots,c_d\in \Q$ with $(q-1)c_i\in \Z$ for all $i$,
 $c_0+c_1+\cdots+c_d\in \Z$ and $c_i\neq 0$ for $i=1,2,\dots,d$.
 Then we define
$$
\widetilde{S}_q(c_0;c_1,\dots,c_d)=
\{c_0+c_1q^{k_1}+c_2q^{k_2}+\cdots+c_dq^{k_d}\mid k_1,k_2,\dots,k_d\in
\N\}
$$
and 
$$S_q(c_0;c_1,\dots,c_d)=\widetilde{S}_q(c_0;c_1,\dots,c_d)\cap
\N.$$
\end{definition}
The conditions on $c_0,c_1,\dots,c_d$ imply that
$$
(q-1)(c_0+c_1q^{k_1}+c_2q^{k_2}+\cdots+c_dq^{k_d})=
$$
$$
=
(q-1)c_0+((q-1)c_1)q^{k_1}+((q-1)c_2)q^{k_2}+\cdots+((q-1)c_d)q^{k_d}\equiv
$$
$$
\equiv
(q-1)c_0+(q-1)c_1+\cdots+(q-1)c_d=
(q-1)(c_0+c_1+\cdots+c_d)\equiv 0 \bmod (q-1).
$$
It follows that $\widetilde{S}_q(c_0;c_1,\dots,c_d)\subseteq \Z$.
If $c_1,c_2,\dots,c_d$ are all negative, then $S_q(c_0;c_1,\dots,c_d)$ is
finite.
 \begin{definition}
 If $c_i>0$ for some $i$ with $1\leq i\leq d$, then 
 $S_q(c_0;c_1,\dots,c_d)$
 is called an {\it elementary $p$-nested set of order $d$}.
A {\it $p$-nested set of order $\leq d$} is 
a union of a finite set and finitely many elementary $p$-nested
sets of order $\leq d$. A $p$-nested set of order $\leq d$ is
said to have order $d$ if is {\it not\/} a $p$-nested set of
order $\leq d-1$.
  \end{definition}

We say that two sets $S,T$ are equal {\it up to a finite set\/} if 
the symmetric difference $(S\setminus T)\cup (T\setminus S)$ is finite.
\begin{definition}\label{defpnormal}
We will call a subset of $\N$ {\it $p$-normal of order $d$} if it is,
up to a finite set, equal to the union of a $p$-nested
set of order $d$ and finitely many infinite arithmetic progressions.
\end{definition}
We are now ready to state the main results of this paper.
\begin{theorem}\label{thm1a}
Suppose that $K$ is a field of characteristic $p>0$. If
$a\in K^{\N}$ is a $K$-recurrence sequence of order $d$, then
${\mathcal Z}(a)$ is $p$-normal of order $\leq d-2$.
\end{theorem}
In the proof
we will first show  that ${\mathcal Z}(a)$ is a $p$-automatic set (a useful
notion from theoretical computer science),
using a technique reminiscent of Frobenius splitting.
The structure of the automaton that produces ${\mathcal Z}(a)$
turns out to be very special. Using this
we will be able to show that ${\mathcal Z}(a)$ is $p$-normal.
Our approach in positive characteristic is entirely different
from the techniques in the proof of the Skolem-Mahler-Lech
theorem in characteristic 0. 

Theorem~\ref{thm1a} can be generalized to recurrence sequences
in modules over $\F_p$-algebras.
\begin{theorem}\label{thm1b}
Suppose that $p$ is a prime number, $R$ is an $\F_p$-algebra, 
$M$ is a left $R$-module
and $a\in M^{\N}$ is an $R$-recurrence sequence of order $d$.
 Then ${\mathcal Z}(a)$ is $p$-normal of order $\leq d-2$.
 \end{theorem}
We will reduce Theorem~\ref{thm1b} to
 Theorem~\ref{thm1a}.

An advantage of our proof of Theorem~\ref{thm1a}
is that, unlike in characteristic 0, all the steps in the proof are effective:
\begin{theorem}\label{theoremeffective}
Let $p$ be a prime. Given a field  $K$  which is finitely generated over $\F_p$ 
and a $K$-recurrence sequence
$a\in K^{\N}$, we can effectively determine ${\mathcal Z}(a)$.
\end{theorem}
In other words there exists an algorithm which,
given $K$ and an explicit
recurrence relation for the sequence $a\in K^{\N}$, 
produces ${\mathcal Z}(a)$ in finite time. The
format of the output is an explicit description of ${\mathcal Z}(a)$
in terms of finite sets, arithmetic progressions and elementary $p$-nested
sets as in the definition of a $p$-normal set (see Definition~\ref{defpnormal}).

Some results are known about the density of the
zeroes of recurrence sequences in positive characteristic.
For a subset $S\subseteq \N$, define
$$
\delta_S(n):=\max_{m\in \N}|S\cap \{m,m+1,\dots,m+n-1\}|.
$$
The {\it upper Banach density} of $S$ is defined by
$$
\delta^+(S)=\limsup_{n\to \infty}\frac{\delta(n)}{n}.
$$
Suppose that $K$ is a field of characteristic $p>0$ and
$a\in K^{\N}$ is a $K$-recurrence sequence of order $d$
such that ${\mathcal Z}(a)$ does not contain any infinite arithmetic
progressions. It was proved in \cite{Bez} that $\delta^+({\mathcal Z}(a))=0$.
In other words, we have $\delta_{{\mathcal Z}(a)}(n)=o(n)$.
(We use here the standard 
{\it big $O$, little $o$, $\Omega$, $\omega$} notations:
$$
\begin{array}{rcl}
f(n)=O(g(n)) & \mbox{ means } & \limsup_{n\to
\infty}\left|f(n)/g(n)\right|<\infty,\\
f(n)=o(g(n)) & \mbox{ means } & \lim_{n\to
\infty}f(n)/g(n)=0,\\
f(n)=\Omega(g(n)) & \mbox{ means } & \liminf_{n\to
\infty}\left|f(n)/g(n)\right|>0,\\
f(n)=\omega(g(n)) & \mbox{ means } & \lim_{n\to
\infty} g(n)/f(n)=0
\end{array}.
$$
for any real-valued functions $f,g\in \R^{\N}$.)

 Indeed, if $\delta^+({\mathcal Z}(a))>0$,
then by Szemer\'edi's theorem (see~\cite{Sz,Fu}) ${\mathcal Z}(a)$  contains
arithmetic progressions of arbitrary length. But then ${\mathcal Z}(a)$
 contains an infinite arithmetic progression (see
Corollary~\ref{corar}),
which contradicts our assumptions.
From effective estimates for $\delta_{{\mathcal Z}(a)}(n)$ in  
\cite{SS1,SS2} and \cite[Theorem 5.9]{GAIT} follows
that 
$$\delta_{{\mathcal Z}(a)}(n)=O(n^{1-\Delta(n)}),
$$
where $\Delta(d)$  is defined by $\Delta(2)=1$ and
$$
\Delta(d+1)=\frac{\Delta(\lfloor\frac{d}{2}\rfloor+1)}{\lfloor\frac{d}{2}\rfloor+1}$$ 
for $d\geq 2$.
One can show that $\Delta(d)\geq \exp(-c\log^2(d))$ for some
constant $c$ (see~\cite{SS1}). In \cite[\S2.5]{GAIT} it was suggested that $\delta_{{\mathcal
Z}(a)}(n)$ might have a logarithmic
upper bound. Although an upper bound
 $O(\log(n))$ is impossible because of Example~\ref{ex2a}, 
Theorem~\ref{thm1a} implies the following result.
\begin{corollary}
Suppose that $K$ is a field of characteristic $p>0$
and $a\in K^{\N}$ is an $K$-recurrence sequence of order $d$.
If ${\mathcal Z}(a)$ does not contain any infinite arithmetic
progressions, then
$$
\delta_{{\mathcal Z}(a)}(n)=O(\log(n)^{d-2}).
$$
\end{corollary}
{\bf Acknowledgement.} I would like to thank Jason Bell for inspiring
discussions, Andreas Blass for help on finite automata and
regular languages, and Jeff Lagarias and David Masser for useful comments.

\section{Preliminaries}
In this section we give some definitions and elementary
facts about linear recurrence sequences. 
Let $R$ be a ring and $M$ be an $R$-module.
We define the {\it shift operator\/} $E:M^{\N}\to M^{\N}$ by 
$$
(Ea)(n):=a(n+1),\quad n\in \N,
$$
for all $a\in M^{\N}$. Scalar multiplication makes $M^{\N}$ into
a left $R$-module. Using the shift operator, we may view $M^{\N}$
as a left $R[E]$-module, where $R[E]$ is the polynomial ring over $R$.
Suppose that $a\in M^{\N}$. Then
the recurrence relation~(\ref{eqrec}) is equivalent to
$$
P(E)\cdot a=0,
$$
where
$$
P(E)=E^d-\gamma_{d-1}E^{d-1}-\gamma_{d-2}E^{d-2}-\cdots-\gamma_1E-\gamma_0.
$$
We call $P(E)$ the {\it companion polynomial} of the recurrence
relation~(\ref{eqrec}).

For $j,k\in \N$ we define $L^k_j:\N\to \N$ by 
$$L^k_j(n)=kn+j.$$
We define the operator $T^k_j:M^\N\to M^\N$ by 
$$
T^k_ja=a\circ L^k_j.
$$
So we have
$$
(T^k_ja)(n)=a(L^k_j(n))=a(kn+j).
$$
Note that $E=T^1_1$. We  have the  relation
\begin{equation}\label{eqcomm}
ET^m_j=T^m_jE^m.
\end{equation}

Suppose that $R$ is a commutative ring and $P(T),Q(T)\in R[T]$ are given by
$$
P(T)=\alpha_nT^n+\alpha_{n-1}T^{n-1}+\cdots+\alpha_0
$$
$$
Q(T)=\beta_mT^m+\beta_{m-1}T^{m-1}+\cdots+\beta_0
$$
where $\alpha_0,\dots,\alpha_n,\beta_0,\dots,\beta_m\in R$ and
$\alpha_n,\beta_m\neq 0$. The {\it resultant} $\Res_T(P(T),Q(T))$ is defined
as the determinant of the matrix
$$
\begin{pmatrix}
\alpha_n & \alpha_{n-1} & \cdots & \alpha_0 & & &\\
	& \alpha_{n} & \alpha_{n-1}& \cdots & \alpha_0 & &\\
	&             & \ddots & &  & \ddots &\\
	&              &       & \alpha_n & \alpha_{n-1} & \cdots & \alpha_0\\
\beta_m & \beta_{m-1} & \cdots & \beta_0 & & &\\
	& \beta_{m} & \beta_{m-1}& \cdots & \beta_0 & &\\
	&             & \ddots & &  & \ddots &\\
	&              &       & \beta_n & \beta_{n-1} & \cdots & \beta_0\\

\end{pmatrix}
$$
(see for example~\cite[IV,\S8]{Lang}).
\begin{lemma}\label{lem2.1}
If $a\in M^{\N}$ is an $R$-recurrence sequence of order $d$, then
$T^k_ja$ is an $R$-recurrence sequence of order $\leq d$.
\end{lemma}
\begin{proof}
See~\cite[Theorem 1.3]{GAIT} in the special case where $R$ is a field. 
Suppose that $P(E)\in R[E]$ is a monic
polynomial of degree $d$ such that $P(E)\cdot a=0$. 
Define $Q(E),U(E)\in R[E]$
by
$$
Q(E)=\Res_F(F^k-E,P(F))
$$
and
$$
U(E)=\Res_F(F^{k-1}+F^{k-2}E+\cdots+E^{k-1},P(F)).
$$
The polynomials $(-1)^{d(k+1)}Q(E)$ and $(-1)^{d(k+1)}U(E)$
are monic  and they have degrees $d$ and $d(k-1)$ respectively.
Using
$$
F^k-E^k=(F-E)(F^{k-1}+F^{k-2}E+\cdots+E^{k-1})
$$
and the multiplicative
property of the resultant (\cite[IX, Theorem 3.10]{Lang}), we get
$$
Q(E^k)=\Res_F(F^k-E^k,P(F))=$$
$$
=\Res_F(F-E,P(F))
\Res_F(F^{k-1}+F^{k-2}E+\cdots+E^{k-1},P(F))=P(E)U(E).
$$
From (\ref{eqcomm}) follows that
$$
Q(E)T^k_ja=T^k_jQ(E^k)a=T^k_jU(E)P(E)a=0.
$$
Therefore, $T^k_ja$ satisfies a recurrence relation of order $d$ with
companion polynomial is $(-1)^{d(k+1)}Q(E)$.
\end{proof}
\begin{corollary}\label{corar}
If $a\in M^{\N}$ is an $R$-recurrence sequence of order $d$, and
$$
a(k)=a(k+m)=a(k+2m)=\cdots=a(k+(d-1)m)=0
$$
for some $k\in \N$ and some positive integer $m$, then
$$
a(k+im)=0
$$
for all $i\in \N$.
\end{corollary}
\begin{proof}
The sequence $b=T^m_ka$ is an $R$-recurrence sequence of order $\leq d$
by Lemma~\ref{lem2.1},
and 
$$
b(0)=b(1)=\cdots=b(d-1)=0.
$$
Using induction and 
the recurrence relation for $b$ we get $a(k+im)=b(i)=0$ for all~$i$.
\end{proof}

Let us assume that $M=R=K$ is a field and let $a\in K^{\N}$. The set
$$
\Ann(a)=\{P(E)\in K[E]\mid P(E) a=0\} 
$$
is a principal ideal in $K[E]$. Therefore, the ideal $\Ann(a)$ is generated
by a unique monic polynomial $P_a(E)$. We call $P_a(E)$ the {\it minimum
polynomial\/} of the recurrence sequence $a$. The degree of $P_a(E)$
is exactly the order of the recurrence sequence. 

Suppose that
\begin{equation}\label{eqmin}
P_a(E)=\prod_{i=1}^r (E-\alpha_i)^{m_i}
\end{equation}
where $\alpha_1,\dots,\alpha_r$ are distinct roots
in the algebraic closure $\overline{K}$ of $K$, and
$m_1,\dots,m_r$ are positive integers.
Then $a$ has order $d=\sum_{i=1}^rm_i$.
 It is well known that $a$
has the form
\begin{equation}\label{eqform}
a(n)=\sum_{i=1}^r\sum_{j=0}^{m_i-1}\beta_{i,j}{n\choose j}\alpha_i^n,\quad n\in
\N,
\end{equation}
where $\beta_{i,j}\in \overline{K}$ for all $i,j$ (see for example~\cite[1.1.6]{GAIT}).
Also, any sequence  $a\in K^{\N}$ of the form (\ref{eqform}) satisfies
a recurrence relation of order $d=\sum_{i=1}^r m_i$ and the companion
polynomial of this recurrence relation is given by~(\ref{eqmin}).
\begin{definition}
The recurrence sequence $a$ is called 
\begin{itemize}
\item  {\it basic\/} if $0$ is not
a root of $P_a(E)$;
\item {\it nondegenerate\/} if all roots of $P_a(E)$ are
nonzero and the quotient of any two distinct roots 
is not a root of unity;
\item {\it simple\/} if all roots of $P_a(E)$ are distinct.
\end{itemize}
\end{definition}
If $a$ is basic, then using the recurrence relation~(\ref{eqrec}) backwards
one can define $a(n)$ for all $n\in \Z$. In that case (\ref{eqform}) would
be valid for all $n\in \Z$ if we interpret
$${n\choose j}=\frac{n(n-1)\cdots (n-j+1)}{j!}
$$
as a polynomial of degree $j$ that is defined for all integers $n$.

If $a$ is simple, then (\ref{eqform}) takes a simpler form, namely
$$
a(n)=\sum_{i=1}^d\beta_i\alpha_i^n,\quad n\in \N,
$$
where $\beta_1,\dots,\beta_d\in \overline{K}$. 
\begin{lemma}\label{lem3.4}
Suppose that $a\in K^{\N}$ is a $K$-recurrence sequence where $K$ is a field.
\begin{enumerate}
\letters
\item There exists an $i$ such that $E^ia$ is basic.
\item If $a$ is basic, then there exists a positive integer $k$ such
that $T^k_{j}a$ is nondegenerate for all $j$.
\item If $a$ is basic and $K$ has characteristic $p>0$,
then there exists a positive integer $k$ such that $T^k_ja$
is simple and nondegenerate for all $j$.
\end{enumerate}
\end{lemma}
\begin{proof} Let $P_a(E)$ be the minimum polynomial
of $a$.

(a)  We can
write $P_a(E)=E^iQ(E)$ where $Q(0)\neq 0$. Then $Q(E)(E^ia)=0$,
hence $E^ia$ is basic.

(b) Let $P_a(E)$ as in (\ref{eqmin}). 
There exists a $k$ such that for all $i\neq j$
 we have that $(\alpha_i/\alpha_j)^k=1$ if and only
if $\alpha_i/\alpha_j$ is a root of unity. 
Define $Q(E)=\Res_F(E-F^k,P(F))$ as in the
proof of Lemma~\ref{lem2.1}. We have $Q(E)(T^k_ja)=0$ for all $j$.
The roots of $Q(E)$ are $\alpha_1^k,\dots,\alpha_r^k$.
In particular, the quotient of any two distinct roots of
$Q(E)$ is not a root of unity different from 1. This shows
that $T^k_ja$ is nondegenerate for all $j$.

(c) Let $\gamma_1,\dots,\gamma_s$ be distinct such that
$$
\{\gamma_1,\dots,\gamma_s\}=\{\alpha_1^k,\dots,\alpha_r^k\}.
$$
Let $q$ be a power of $p$ such that $q\geq d$. Define
$U(E)=\prod_{i=1}^r(E-\gamma_i^q)$. Then $Q(E)$ divides
$$U(E^q)=\prod_{i-1}^r(E^q-\gamma_i^q)=\prod_{i=1}^r(E-\gamma_i)^q.$$
Because $P_a(E)$ divides $Q(E^k)$, it divides $U(E^{qk})$ as well.
Therefore, for all $j$ we have
$$
U(E) (T^{qk}_ja)=T^{qk}_j(U(E^{qk})a)=0.
$$
Note that $U(E)$ has distinct roots, and any quotient
of two distinct roots is not a root of unity. It follows
that $T^{qk}_ja$ is simple and nondegenerate for all $j$.
\end{proof}

\begin{definition}
A {\it $d$-balanced\/} subset of $\N$ is a set of the form
$$
\{m_0+k_1m_1+\cdots+k_dm_d\mid k_1,\dots,k_d\in \{0,1\}\},
$$
where $m_0\in \N$ and $m_1,\dots,m_d$ are positive integers.
\end{definition}
If $m_1=m_2=\cdots=m_d$, then
$$
\{m_0+k_1m_1+\cdots+k_dm_d\mid k_1,\dots,k_d\in
\{0,1\}\}=\{m_0,m_0+m_1,m_0+2m_1,\dots,m_0+dm_1\},
$$
so an arithmetic progression of length $d+1$ is $d$-balanced.
\begin{lemma}\label{lem2.6}
Suppose that $K$ is a field and $a\in K^{\N}$ is a nonzero $K$-recurrence sequence
of order $d$.
\begin{enumerate}
\letters
\item
If $K$ has characteristic 0 and 
$a\in K^{\N}$ is  nondegenerate 
then ${\mathcal Z}(a)$ does not contain a $(d-1)$-balanced subset.
\item
If $K$ has positive characteristic and $a\in K^{\N}$ is simple and
nondegenerate, then ${\mathcal Z}(a)$ does not contain a $(d-1)$-balanced
subset.
\end{enumerate}
\end{lemma}
\begin{proof}
Suppose that
$$
S:=\{m_0+k_1m_1+\cdots+k_dm_d\mid k_1,\dots,k_d\in
\{0,1\}\}\subseteq {\mathcal Z}(a).
$$
Let
$$
P_a(E)=\prod_{i=1}^d(E-\alpha_i)
$$
be the minimum polynomial of $a$.
Define 
\begin{equation}\label{eqQ}
Q(E)=E^{m_0}\prod_{i=1}^{d-1}(E^{m_i}-\alpha_i^{m_i}).
\end{equation}
If $i<d$, then $E-\alpha_d$ does not divide
$$
\frac{E^{m_i}-\alpha_i^{m_i}}{E-\alpha_i}=E^{m_i-1}+\alpha_iE^{m_i-2}+\cdots+\alpha_i^{m_i-1}:
$$
In case (b), $\alpha_d^{m_i}\neq \alpha_i^{m_i}$ since $\alpha_i/\alpha_d$
is not a root of unity for $i=1,2,\dots,d-1$. In case (a), $\alpha_d^{m_i}=\alpha_i^{m_i}$
implies that $\alpha_d=\alpha_i$. Since $K$ has characteristic $0$,
$E^{m_i}-\alpha_i^{m_i}$ is square-free, and $E-\alpha_i=E-\alpha_d$ does
not divide $(E^{m_i}-\alpha_i^{m_i})/(E-\alpha_i)$.

The polynomial $P_a(E)$ divides
$Q(E)(E-\alpha_d)$, but not $Q(E)$.
From $(E-\alpha_d)Q(E) a=0$ follows that
$$
(Q(E) a)(n)=\beta\alpha_d^n
$$
for some $\beta\in \overline{K}$. If the coefficient of $E^i$ in $Q(E)$ is 
nonzero, then (\ref{eqQ}) implies that 
$i\in S$ and $a(i)=0$.
It follows that $(Q(E)a)(0)$ is a linear combination
of $a(i),i\in S$. Therefore $\beta=(Q(E) a)(0)=0$.
But then $Q(E) a=0$ and $Q(E)$ must be divisible by the 
minimum polynomial $P_a(E)$. Contradiction.
\end{proof}
Theorem~\ref{thm1a} can be reduced to the following theorem for
simple nondegenerate sequences.
\begin{theorem}\label{thm1c}
If $K$ is a field of characteristic $p>0$ and $a\in K^{\N}$
is a nonzero simple nondegenerate $K$-recurrence sequence of order $d$, 
then ${\mathcal Z}(a)$ is a $p$-nested set of order $\leq d-2$.
\end{theorem}

\begin{proof}[Proof of Theorem~\ref{thm1a}]
Suppose that $a\in K^{\N}$ is an $K$-recurrence sequence of order $\leq d$.
Without loss of generality we may assume that $a$ is basic, because
any linear recurrence sequence can be changed into a  
basic one by changing only finitely many entries in the sequence.
There exist a positive integer $k$ such that $T^k_ja$ is simple
and nondegenerate of order $\leq d$ for $j=0,1,\dots,k-1$
by Lemma~\ref{lem3.4} and
Lemma~\ref{lem2.1}. We have
\begin{equation}\label{eqza}
{\mathcal Z}(a)=\bigcup_{j=0}^{k-1}L^k_j({\mathcal Z}(T^k_j(a)).
\end{equation}
If $T^k_ja=0$, then ${\mathcal Z}(T^k_ja)=\N$ and
$L^k_j({\mathcal Z}(T^k_ja))=L^k_j(\N)=j+k\N$ is an infinite arithmetic
progression. Otherwise, ${\mathcal Z}(T^k_ja)$ is a
$p$-nested set of order $\leq d-2$ by Theorem~\ref{thm1c}.
But then $L^k_j({\mathcal Z}(T^k_ja))$ is   a
union of a finite set and finitely many sets of the form
$$
L^k_j(S_q(c_0;c_1,\dots,c_r))=S_q(j+kc_0;kc_1,\dots,kc_r).
$$
 with $q$ a power of $p$
and $r\leq d-2$, hence $L^k_j({\mathcal Z}(T^k_ja))$ is
$p$-nested of order $\leq d-2$ for all $j$.
From (\ref{eqza}) follows that ${\mathcal Z}(a)$ is $p$-normal
of order $\leq d-2$ as well.
\end{proof}
\section{Examples in positive characteristic}
In this section we will concentrate on simple nondegenerate
$K$-recurrence sequences in $K$ where $K$ is a field of positive
characteristic. 
The main idea behind the construction of various
pathological examples is the following proposition.
\begin{proposition}\label{proptensor}
Assume that $K$ is a field of characteristic $p>0$ and $q$ is
a power of $p$.
Suppose that
$a\in K^{\N}$ is given by
$$
a(n)=\sum_{i=1}^{d}\beta_i\alpha_i^n,\quad n\in\N,
$$
where $\alpha_1,\dots,\alpha_d,\beta_1,\dots,\beta_d\in \overline{K}$.
If for some $c_0,c_1,\dots,c_r\in \Z$ the sum
$$
\sum_{i=1}^d (\beta_i\alpha_{i}^{c_0})\otimes \alpha_i^{c_1}\otimes \cdots
\otimes \alpha_{i}^{c_r}\in \underbrace{\overline{K}\otimes_{\F_q}\overline{K}\otimes_{\F_q}\otimes \cdots
\otimes_{\F_q}\overline{K}}_{r+1},
$$
is equal to $0$, then 
$$
S_q(c_0;c_1,\dots,c_r)\subseteq {\mathcal Z}(a).
$$
\end{proposition}
\begin{proof}
Let $\phi:\overline{K}\to \overline{K}$ be the Frobenius map defined by
$\phi(\alpha)=\alpha^p$. Define an $\F_q$-linear map
$$\psi:\overline{K}\otimes_{\F_q}\overline{K}\otimes_{\F_q}\cdots
\otimes_{\F_q}\overline{K}\to \overline{K}
$$
by
$$
\gamma_0\otimes \gamma_1\otimes \cdots \otimes \gamma_r\mapsto
\gamma_0\gamma_1\cdots \gamma_r.
$$
We have
$$
a(c_0+c_1q^{k_1}+\cdots+c_rq^{k_r})=\sum_{i=1}^d
\beta_i\alpha_i^{c_0}\alpha_i^{c_1q^{k_1}}\alpha_i^{c_2q^{k_2}}\cdots
\alpha_i^{c_rq^{k_r}}=
$$
$$
=\sum_{i=1}^r\psi\big((\beta_i\alpha_i^{c_0})\otimes\alpha_i^{c_1q^{k_1}}\otimes
\alpha_i^{c_2q^{k_2}}\otimes\cdots
\otimes\alpha_i^{c_rq^{k_r}}\big)=
$$
$$
=\sum_{i=1}^r\psi\big((\beta_i\alpha_i^{c_0})\otimes
\phi^{k_1}(\alpha_i^{c_1})\otimes
\phi^{k_2}(\alpha_i^{c_2})\otimes\cdots
\otimes\phi^{k_r}(\alpha_i^{c_r})\big)=
$$
$$
=\psi\big(\id\otimes \phi^{k_1}\otimes\phi^{k_2}\otimes\cdots\otimes
\phi^{k_r}\big(\sum_{i=1}^r(\beta_i\alpha_{i}^{c_0})\otimes \alpha_i^{c_1}\otimes \cdots
\otimes \alpha_{i}^{c_r}\big)\big)=
$$
$$
=\psi\big(\id\otimes \phi^{k_1}\otimes\phi^{k_2}\otimes\cdots\otimes
\phi^{k_r}(0)\big)=\psi(0)=0.
$$
\end{proof}
\begin{proposition}\label{prop1}
For the sequence $a\in \F_p(x,y,z)$ defined by 
$$
a(n)=(x+y+z)^n-(x+y)^n-(x+z)^n-(y+z)^n+x^n+y^n+z^n,
$$
we have
$$
{\mathcal Z}(a)=\{p^n\mid n\in \N\}\cup \{p^n+p^m\mid n,m\in \N\}.
$$
\end{proposition}
\begin{proof}
From
$$
1\otimes (x+y+z)+(-1)\otimes (x+y)+(-1)\otimes (x+z)+(-1)\otimes (y+z)+
1\otimes x+1\otimes y+1\otimes z=0
$$
in $\F_p(x,y,z)\otimes_{\F_p}\F_p(x,y,z)$ and Proposition~\ref{proptensor}
follows that 
$$
S_p(0;1)=\{1,p,p^2,\dots\}\subseteq {\mathcal Z}(a).
$$
One can also check that
\begin{multline}\label{eq1}
1\otimes (x+y+z)\otimes (x+y+z)+(-1)\otimes (x+y)\otimes (x+y)+
(-1)\otimes (x+z)\otimes (x+z)+\\
(-1)\otimes (y+z)\otimes (y+z)+
1\otimes x\otimes x+1\otimes y\otimes y+1\otimes z\otimes z=0
\end{multline}
in $\F_p(x,y,z)\otimes_{\F_p}\F_p(x,y,z)\otimes_{\F_p}\F_p(x,y,z)$.
Again from Proposition~\ref{proptensor}
follows that 
$$
S_p(0;1,1)=\{p^n+p^m\mid n,m\in \N\}\subseteq {\mathcal Z}(a).
$$

Conversely, if $n$ is not a $p$-power or the sum of 2 powers of $p$
then we can write $n=s+t+ul$ where
$s=p^i$, $t=p^j$, $u=p^k$, $l$ is a positive integer not divisible by $p$
and $i<j<k$.
We get
\begin{multline}
(x+y+z)^n=(x+y+z)^{s+t+ul}=(x+y+z)^s(x+y+z)^t((x+y+z)^u)^l=\\
=(x^s+y^s+z^s)(x^t+y^t+z^t)(x^u+y^u+z^u)^l
\end{multline}
From this it is clear the the coefficient of $x^sy^tz^{ul}$
in $(x+y+z)^n$ is equal to 1, so $a(n)\neq 0$.\end{proof}
The previous example easily generalizes to the following one.
\begin{example}\label{ex5}
Define $a_d\in \F_p(x_1,\dots,x_d)^{\N}$ by
$$
a_d(n):=\sum_{I}
(-1)^{d+|I|}\big(\sum_{i\in I}x_i\big)^n
$$
where $I$ runs over all nonempty subsets of $\{1,2,\dots,d\}$ and
$|I|$ denotes the cardinality of $I$.
The set ${\mathcal Z}(a_d)$ consists of all sums
of at most $d$ powers of $p$.
\end{example}
The phenomenon of Proposition~\ref{prop1} already appears in recurrence sequences
of order 4 as the following example shows.
\begin{example}\label{ex6}
Consider the sequence $a\in \F_4(x)^{\N}$ defined by
$$
a(n)=x^n+(x+1)^n+(x+\alpha)^n+(x+1+\alpha)^n,
$$
where $\alpha\in \F_4\setminus\F_2$. We can compute
$a(0)=a(1)=a(2)=0$ and $a(3)=1$.
This  sequence satisfies a recurrence relation of order 4, whose
companion polynomial is
$$
(E-x)(E-x-1)(E-x-\alpha)(E-x-1-\alpha)=E^4+E+(x^4+x).
$$
The  recurrence relation for $a$ is
$$
a(n+4)=a(n+1)+(x^4+x)a(n).
$$
Note that $a$ is actually a sequence in the subfield $\F_2(x)\subseteq \F_4(x)$.

We have
$$
{\mathcal Z}(a)=\{4^n+4^m\mid n,m\in \N\}\cup \{2\cdot (4^n+4^m)\mid n,m\in
\N\}\cup \{0,1\}.
$$
\end{example}
In the remainder of this section we describe a construction of
simple nondegenerate recurrence sequences with many zeroes.
Suppose that $p$ is a prime and $q$ is a power of $p$.
For $c_0,c_1,\dots,c_r\in \Z$,  we will define
a nonzero simple nondegenerate sequence $a\in \F_q(x)^{\N}$
such that
$$
S_q(c_0;c_1,\dots,c_r)\subseteq {\mathcal Z}(a)
$$
Let $\gamma_1(x),\gamma_2(x),\dots,\gamma_{s(d)}(x)$ be all the 
irreducible polynomials
in $\F_q[x]$ of degree $\leq d$. Define
$$
{\mathcal C}_d:=\big\{\textstyle\prod_{i=1}^{s(d)}\gamma_i(x)^{e_i}\big| e_i
\in \{0,1\}\mbox{ for all
$i$}\big\}.
$$
The set ${\mathcal C}_d$ has $2^{s(d)}$ elements.
Let $\F_q[x]_{\leq e}$ be the space of polynomials of degree $\leq e$.
Define ${\mathcal M}_{d,c}$ as the $\F_q$-vector space
spanned by all $\gamma(x)^c$ with $\gamma(x)\in {\mathcal C}_d$.
If $c$ is a nonnegative integer, then ${\mathcal M}_{d,c}$ is contained
in $\F_q[x]_{\leq cds(d)}$ and therefore
$$
\dim {\mathcal M}_{d,c}\leq \dim \F_q[x]_{\leq cds(d)}=cds(d)+1.
$$
If $c<0$ then ${\mathcal M}_{d,c}$ is contained
in ${\mathcal M}_{d,-c}(\prod_{i=1}^{s(d)}\gamma_i(x)^{c})$.
So we have
$$
\dim {\mathcal M}_{d,c}\leq |c|ds(d)+1
$$
for all integers $c$. For all $c\in \Z$, the inequality
$$
\dim \F_q[x]_{\leq e}{\mathcal M}_{d,c}\leq |c|ds(d)+e+1
$$
holds. Consider the vector space
$$
V_{d,e}=\F_q[x]_{\leq e}{\mathcal M}_{d,c_0}\otimes_{\F_q} {\mathcal
M}_{d,c_1}\otimes_{\F_q}{\mathcal M}_{d,c_2}\otimes_{\F_q} \cdots
\otimes_{\F_q}{\mathcal M}_{d,c_r}.
$$
We have
$$
\dim V_{d,e}\leq (|c_0|ds(d)+e+1) \prod_{i=1}^r (|c_i|ds(d)+1).
$$
Let ${\mathcal B}_{d,e}\subseteq V_{d,e}$ be the subset of all
$$
(x^i\gamma(x)^{c_0})\otimes \gamma(x)^{c_1}\otimes \cdots \otimes \gamma(x)^{c_r}
$$
with $0\leq i\leq e$ and $\gamma(x)\in {\mathcal C}_d$.
The set ${\mathcal B}_{d,e}$ has cardinality $|{\mathcal B}_{d,e}|=(e+1)2^{s(d)}$.

Note that $\lim_{d\to \infty}s(d)=\infty$.
  We have
$$
|{\mathcal B}_{d,e}|>\dim V_{d,e}
$$
for $d\gg 0$
because $|{\mathcal B}_{d,e}|$
depends at least exponentially on $s(d)$ and
 $\dim V_{d,e}$ depends at most polynomially on $s(d)$.
For $d$ large enough, ${\mathcal B}_{d,e}$ will be dependent.
This means that there exist polynomials $v_{\gamma}(x)\in \F_q[x]$,
of degree $\leq e$, not all 0, such that
$$
\sum_{\gamma(x)\in {\mathcal C}_{d}}
(v_{\gamma}(x)\gamma(x)^{c_0})\otimes \gamma(x)^{c_1}\otimes \cdots \otimes 
\gamma(x)^{c_r}=0.
$$
From Proposition~\ref{proptensor} follows that
$$
S_q(c_0;c_1,\dots,c_r)\subseteq {\mathcal Z}(a)
$$
where $a\in \F_q[x]^{\N}$ is defined by
$$
a(n)=\sum_{\gamma(x)\in {\mathcal C}_d}v_{\gamma(x)} \gamma(x)^{n}.
$$
\begin{conjecture}\label{conj1}
If we choose $e$ and $d$ large enough, then there exists
a choice for the $\{v_{\gamma(x)}\}$ such that
$$
S_q(c_0;c_1,\dots,c_r)={\mathcal Z}(a).
$$
\end{conjecture}
Perhaps instead of choosing the $v_{\gamma(x)}\in \F_q[x]_{\leq e}$
one should choose them in $K[x]_{\leq e}$ where $K$ is a field
extension of $\F_q$ containing many transcendental elements.

The converse of Theorem~\ref{thm1a} is not  true. For example, 
consider the $2$-normal set
$$
S=\{2^n\mid n\in \N\}\cup \{2^n+2^m\mid n,m\in \N\}=\{1\}\cup \{2^n+2^m\mid
n,m\in \N\}.
$$
The set $S$ is $2$-normal of order $2$. However, there does
not exist a recurrence sequence $a\in K^{\N}$ with ${\mathcal Z}(a)=S$
of order $\leq 4$: Note that $S$ contains $\{1,2,3,4,5,6\}$
which is an arithmetic
progression of length $6$. Suppose that $a\in K^{\N}$ is a $K$-recurrence
sequence
where $K$ is a field of characteristic $0$.
If  ${\mathcal Z}(a)=S$, then ${\mathcal Z}(a)$ does
not contain an infinite arithmetic progression and $a$ has
order $\geq 7$ by Corollary~\ref{corar}.
However, we do conjecture the following weaker converse.
\begin{conjecture}\label{conj2}
If $S$ is a $p$-normal set, then there exists a 
field $K$ of characteristic $p$ and a $K$-recurrence sequence
$a\in K^{\N}$ such that ${\mathcal Z}(a)=S$.
\end{conjecture}
\begin{lemma}
Conjecture~\ref{conj1} implies Conjecture~\ref{conj2}.
\end{lemma}
\begin{proof}
Suppose Conjecture~\ref{conj1} is true. We will prove Conjecture~\ref{conj2}.
If we change finitely many entries in  a recurrence sequence then
the sequence remains a recurrence sequence.
It follows that if $S,T$ are equal up to a finite set and $T$ is the
zero set of a recurrence sequence in $K^\N$, then
so is $S$.  Without loss of generality we may assume
that $S$ is a  union of finitely many infinite arithmetic progressions
and a $p$-nested set.

Note that if $a,b\in K^{\N}$
are $K$-recurrence sequences, then so is the product $ab$ defined by
$ab(n)=a(n)b(n)$. In particular
$$
{\mathcal Z}(ab)={\mathcal Z}(a)\cup {\mathcal Z}(b).
$$
So we easily reduce to the case where $S$ is either
an arithmetic progression, or an elementary $p$-nested set.

For any arithmetic progression it is easy
to find a recurrence sequence with that particular zero set. 
Suppose that $S=S_q(c_0;c_1,\dots,c_r)$.
If $c_0,\dots,c_r$ are all integers then  Conjecture~\ref{conj1}
implies Conjecture~\ref{conj2}.
 Otherwise, we still have that $(q-1)c_i\in \Z$ for all $i$.
 There exists a $K$-recurrence sequence $a\in K^{\N}$
 such that
 $$
 {\mathcal Z}(a)=S_q((q-1)c_0;(q-1)c_1,\dots,(q-1)c_r)
 $$
 by Conjecture~\ref{conj1}.
Define $b\in K^{\N}$ by $b(n)=a((q-1)n)$ for all $n\in \N$, i.e.,
$b=T^{q-1}_0a$. It follows that
$$
{\mathcal Z}(b)=(L^{q-1}_0)^{-1}({\mathcal Z}(a))=$$
$$=
(L^{q-1}_0)^{-1}(S_q((q-1)c_0;(q-1)c_1,\dots,(q-1)c_r))=
S_q(c_0;c_1,\dots,c_r).
$$
\end{proof}
\section{$p$-Automatic sequences}
We first give the necessary definitions for finite automata.
Let $p$ be a positive integer and define the {\it alphabet\/} ${\mathcal
A}=\{0,1,2,\dots,p-1\}$. Let ${\mathcal A}^\star$ be the set 
of all finite words in the alphabet ${\mathcal A}$ (including the empty word).
 A subset ${\mathcal L}\subseteq {\mathcal A}^\star$
is called a {\it language}. 
\begin{definition}\label{defreg}
(See \cite[\S 1.9]{LP}.)
The set of {\it regular} languages ${\mathcal R}$ is the smallest subset of the
set of all languages such that
\begin{enumerate}
\item $\emptyset\in {\mathcal R}$, and $\{v\}\in {\mathcal R}$ for all $v\in
{\mathcal A}$;
\item If $L,M\in {\mathcal R}$, then $L\cup M\in {\mathcal R}$,
$L\circ M\in {\mathcal R}$ (the set of all concatenations 
of a word in $L$ and a word in $M$).
\item If $L\in {\mathcal R}$, then $L^\star\in {\mathcal R}$ where
$L^\star$ is the Kleene closure, i.e., the set of
all possible concatenations of elements in $L$.
\end{enumerate}
\end{definition}
\begin{definition}
A {\it finite automaton\/} with alphabet ${\mathcal A}$
is a finite set ${\mathcal V}$ called the {\it set of states}, an initial state $I\in {\mathcal V}$, a set of
final states ${\mathcal F}\subseteq {\mathcal V}$ together with a map
$$
\tau:{\mathcal V}\times {\mathcal A}\to {\mathcal V}.
$$
\end{definition}
We will write $S \cdot t$ instead of $\tau(S,t)$ for $t\in {\mathcal A}$
and $S\in {\mathcal V}$. For a word $w=t_rt_{r-1}\cdots t_0$ we inductively define

$$
S\cdot w=S\cdot t_rt_{r-1}\cdots t_0:=(S\cdot t_r)\cdot (t_{r-1}t_{r-2}\cdots
t_0).
$$
This way, we may view $\tau$ as a right action of the monoid ${\mathcal A}^\star$
on ${\mathcal V}$.
We say that the automaton {\it accepts\/} the word $w$ if 
$I\cdot w\in {\mathcal F}$.
\begin{theorem} (See~\cite[Theorem~2.5.1]{LP})\label{thm4.2}
A language ${\mathcal L}$ is regular if and only if it is
the set  accepted words for some automaton.
\end{theorem}

An automaton can be represented by a graph. The set
of vertices are labeled by ${\mathcal V}$. For each state $S\in {\mathcal V}$
and each $i\in {\mathcal A}$ we draw an arrow from $S$ to $S\cdot i$
with label $i$. The initial state vertex we will draw as a square and
all other states will be round. Final states are solid vertices and
all other states are open vertices.
 \begin{example}
Consider an automaton for $p=2$ with ${\mathcal V}=\{I,A,B\}$ where $I$ is an
initial state and $A$ and $B$ are two other states. The set of
final states is ${\mathcal F}=\{A\}$. We define $\tau:{\mathcal V}\times
\{0,1\}\to {\mathcal V}$ by
$I\cdot 0=A$,
$I\cdot 1=B$, $A\cdot 0=B$, $A\cdot 1=I$, $B\cdot 0=I$ and $B\cdot 1=A$.
The set of accepted words is
$$
{\mathcal L}=\{0,11,001,010,100,0000,0111,1011,1101,1110,\dots\}.
$$
In fact, ${\mathcal L}$ is the set of all words
for which the number of $0$'s minus the number of $1$'s is congruent
1 modulo $3$. The graph of the automaton is as follows:\\[10pt]
\centerline{\includegraphics[scale=.4]{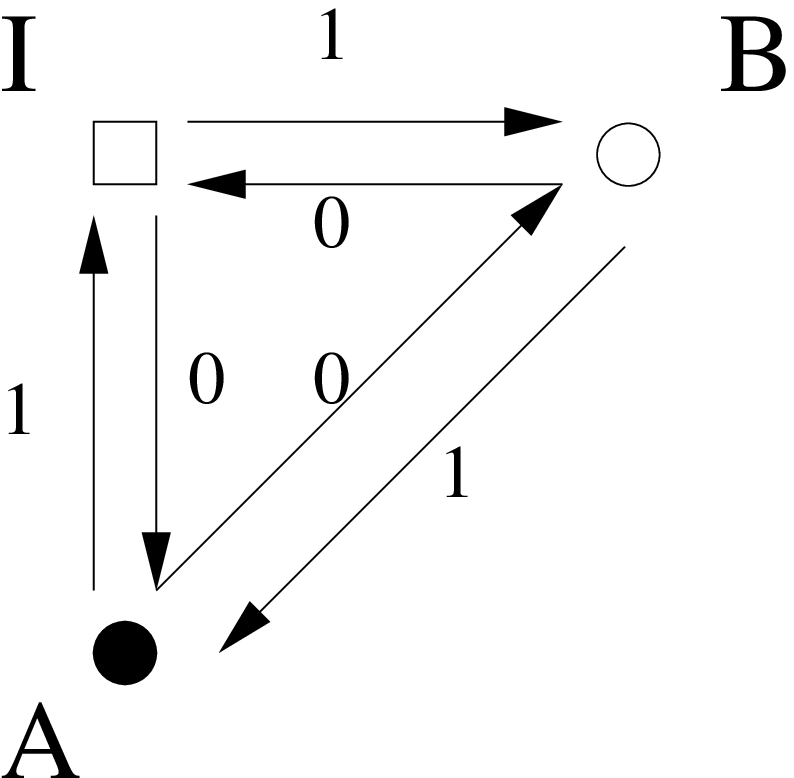}}
\end{example}
Words in the alphabet ${\mathcal A}$ can be viewed as nonnegative
integers written in base $p$.
 For a word $w=t_rt_{r-1}\cdots t_0$ we define
$$
[w]_p=t_rp^r+t_{r-1}p^{r-1}+\cdots+t_0.
$$
For a nonnegative integer $n\in \N$ there exists a unique word $w=t_rt_{r-1}\cdots
t_0\in {\mathcal A}^\star$ with $t_r\neq 0$ and $n=[w]_p$.
\begin{definition}
A subset of $S\subseteq \N$ is called {\it $p$-automatic\/} if there exists
an automaton such that a word $w\in {\mathcal A}^\star$ is accepted if and 
only if $[w]_p\in S$. We call such an automaton an {\it
automaton that produces $S$}.
\end{definition}
In other words, $S\subseteq \N$ is $p$-automatic if and only if the language
$\{w\mid [w]_p\in S\}$ is regular by Theorem~\ref{thm4.2}.
Note that such an automaton accepts a word $w$ if and only
if it accepts the word $0w$. Without loss of generality one
may assume that $I\cdot 0=I$.

If $w=t_rt_{r-1}\cdots t_0$ is a word, then the reverse word $w^{\rm rev}$
is defined by
$$
w^{\rm rev}:=t_0t_1\cdots t_r.
$$
We also define 
$$
w\cdot S=S\cdot w^{\rm rev}
$$
for any word $w\in {\mathcal A}^\star$ and any state $S\in
{\mathcal V}$. This way, we can view ${\mathcal A}^\star$
as a monoid acting on the {\it left\/} on ${\mathcal V}$.
\begin{definition}
A subset of $S\subseteq \N$ is called {\it reversely $p$-automatic\/} 
if there exists
an automaton such that a word $w\in {\mathcal A}^\star$ is accepted if and 
only if $[w^{\rm rev}]_p\in S$. We call such an automaton an {\it
automaton that produces $S$ reversely}.
\end{definition}

\begin{lemma}\label{Blass}
A subset of $S\subseteq \N$ is $p$-automatic if and only if it
is reversely $p$-automatic.
 \end{lemma}
 \begin{proof}
 Let ${\mathcal L}=\{w\in {\mathcal A}^\star\mid [w]_p\in S\}$.
 Now $S$ is $p$-automatic if and only if ${\mathcal L}$ is regular,
 and $S$ is reversely $p$-automatic if and only if 
 $${\mathcal L}^{\rm rev}:=\{w^{\rm rev}\mid w\in {\mathcal L}\}$$
 is regular.
  It is clear from the symmetry in Definition~\ref{defreg} that
 ${\mathcal L}$ is regular if and only if ${\mathcal L}^{\rm rev}$
 is regular. (I thank Andreas Blass for pointing this out to me.)
 \end{proof}
Suppose that a $p$-automaton  produces $S\subseteq \N$. Then the
set $(L^p_t)^{-1}(S)$ is also $p$-automatic. Indeed, if instead of $I$,
we take $I\cdot t$ as initial state, then the automaton
will produce $(L^p_{t})^{-1}(S)$ instead of $S$.
If $w=t_rt_{r-1}\cdots t_0$ is a word, then
$$
(L^p_{t_r})^{-1}(L^{p}_{t_{r-1}})^{-1}\cdots (L^p_{t_0})^{-1}(S)=
(L^p_{t_0}L^p_{t_1}\cdots L^p_{t_r})^{-1}(S)=(L^{p^{r+1}}_{[t_rt_{r-1}\cdots
t_0]_p})^{-1}(S)
$$
is $p$-automatic as well. It is produced by the same automaton
that produces $S$ except that we change the initial state to
$I\cdot t_0t_1\cdots t_r$. Since there only finitely many states,
we get the following corollary.
\begin{corollary}\label{corVS}
If $S\subseteq \N$ is $p$-automatic then the set
$$
{\mathcal V}_S:=\{(L_i^{p^r})^{-1}(S)\mid r\in \N, 0\leq i<p^r\}
$$
is finite.
\end{corollary}
We will also prove the converse.
Suppose that ${\mathcal V}_S$ is finite for some $S\subseteq \N$.
Rather than constructing an automaton that produces $S$, we will construct
an automaton that produces $S$ reversely. In fact, we can do
this in a canonical way. For the set of states we take ${\mathcal V}_S$.
For the initial state we take $I:=S$. For $t\in {\mathcal A}$
and $U\in {\mathcal V}_S$ we define $U\cdot t:=(L^p_t)^{-1}(U)$.
The set of final states is 
$$
{\mathcal F}_S:=\{U\in {\mathcal V}_S\mid 0\in U\}.
$$
A word $w=t_0t_1\dots t_r\in {\mathcal A}^\star$ is accepted
by this automaton if and only if
$$
S\cdot t_0t_1\cdots t_r=(L^p_{t_r})^{-1}(L^{p}_{t_{r-1}})^{-1}\cdots (L^p_{t_0})^{-1}(S)=
(L^{p^{r+1}}_{[t_rt_{r-1}\cdots
t_0]_p})^{-1}(S)=(L^{p^{r+1}}_{[w^{\rm rev}]_p})^{-1}(S)
$$
contains $0$. Therefore $w$ is accepted if and only if $[w^{\rm rev}]_p\in S$.
\begin{proposition}
The converse of Corollary~\ref{corVS} is true: if ${\mathcal V}_S$ is finite,
then $S$ is $p$-automatic. 
The  automaton constructed above has the smallest number of states (namely $|{\mathcal V}_S|$)
among all automata that produce $S$ reversely.
\end{proposition}
\begin{proof}
Suppose that ${\mathcal V}_S$ is finite.
We already constructed an automaton that produces $S$ reversely
whose set of states is ${\mathcal V}_S$.
Assume that we have another automaton that produces $S$ reversely.
Suppose that this automaton is given by
 the set of states ${\mathcal V}$, the initial state $I\in {\mathcal V}$,
 a set of final states ${\mathcal F}\subseteq {\mathcal V}$ and
$$
\tau:{\mathcal V}\times {\mathcal A}\to {\mathcal V}.
$$
By definition $w\cdot I\in {\mathcal F}$ if and only if $[w]_p\in S$.
Without loss of generality we may assume that each state in ${\mathcal V}$
can be reached by a path from $I$.
We define a map $\psi:{\mathcal V}\to {\mathcal V}_S$ as follows. 
For any state $R\in {\mathcal V}$ we define
$$
\psi(R)=\{[w]_p\mid w\cdot R\in {\mathcal F}\}.
$$
If $R=u\cdot I$ and $u$ has length $r$ then
$$
[wu]_p\in S\Leftrightarrow wu\cdot I\in {\mathcal F}\Leftrightarrow
w\cdot R\in {\mathcal F}\Leftrightarrow [w]_p\in \psi(R).
$$
It follows that
$$
\psi(R)=(L_{[u]_p}^{p^r})^{-1}(S)\in {\mathcal V}_S.
$$
From this it is easy to see that $\psi$ is well-defined and surjective.
\end{proof}

\begin{example}
Suppose that 
$$S=\{0,1,3,4,7,9,10,12,16,\dots\}\subseteq \N$$ 
is the subset of all nonnegative integers
which have an even number of $0$'s in their binary expansion, together with 0.
If we write integers in base 2 we have
$$
S=\{0,1,11,100,111,1001,1010,1100,\dots\}.
$$
The set $S$ is  $2$-automatic. An automaton that produces
$S$ is for example:\\[10pt]
\centerline{\includegraphics[scale=.4]{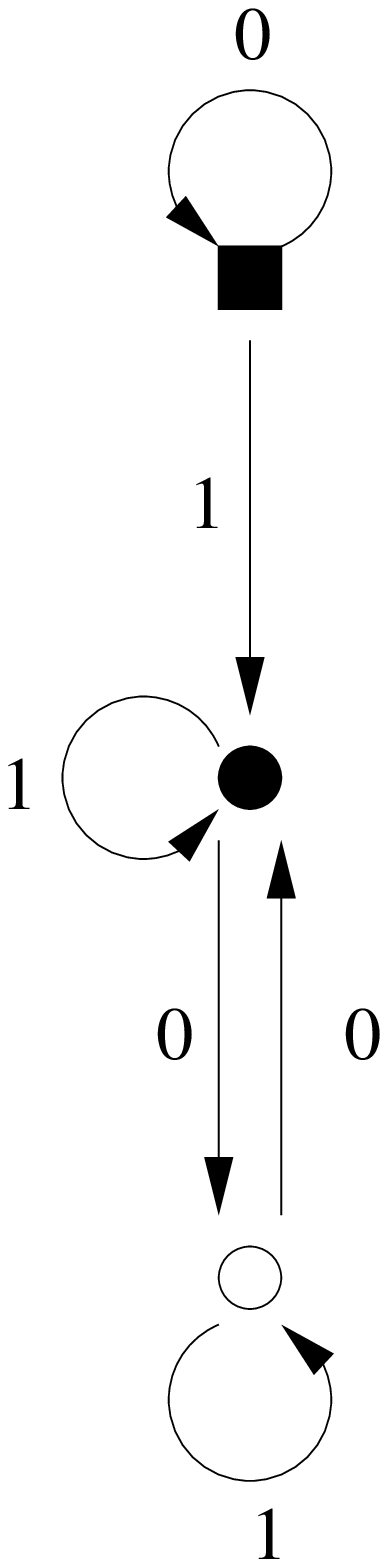}}
The set ${\mathcal V}_S$ contains $4$ elements. The first element is $S$ itself.
The set $(L^{2}_0)^{-1}(S)$ consists of all positive integers
with an odd number of zeroes in their binary expansion together with 0.
The set $(L^4_2)^{-1}(S)$ is the set of all {\it positive\/}
integers with an odd numbers of zeroes in their binary expansion.
 Finally, the set $(L^8_2)^{-1}(S)$ consists of
all positive integers with an even number of zeroes in their
binary expansion. Using the set ${\mathcal V}_S$ we construct
a automaton that produces $S$ reversely. The graph of the automaton
is as follows:\\[10pt]
\centerline{\includegraphics[scale=.4]{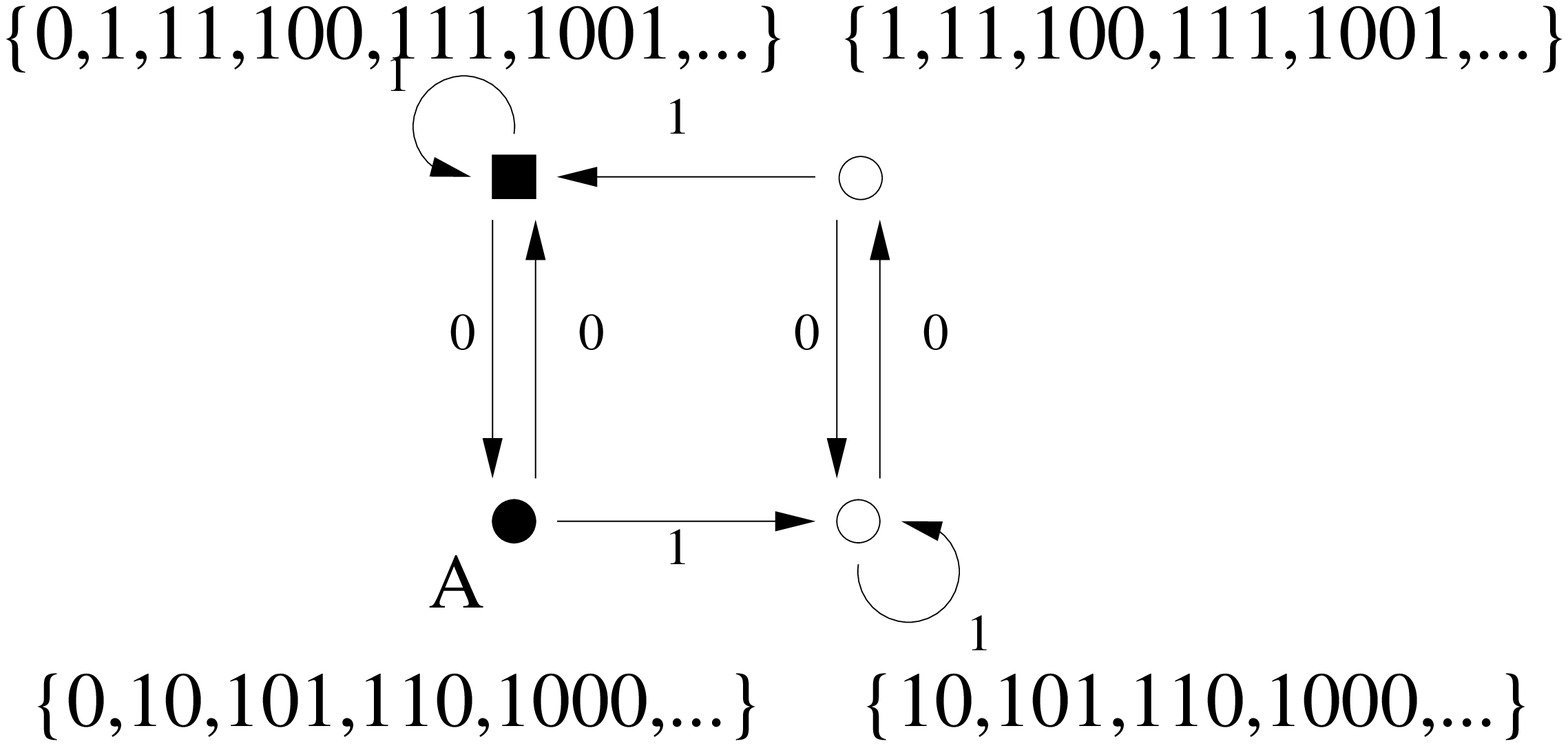}}
\end{example}
\begin{example}
Let $p=2$ and let 
$$S=\{0,3,5,6,9,10,12,15,31,\dots\}$$ 
be the set of all 
nonnegative
integers
which have an even number of $1$'s in their base $2$ expansion.
In base 2 we get
$$
S=\{0,11,101,110,1001,1010,1100,1111,10001,\dots\}.
$$
The  automaton producing $S$ reversely is as follows:\\[10pt]
\centerline{\includegraphics[scale=.4]{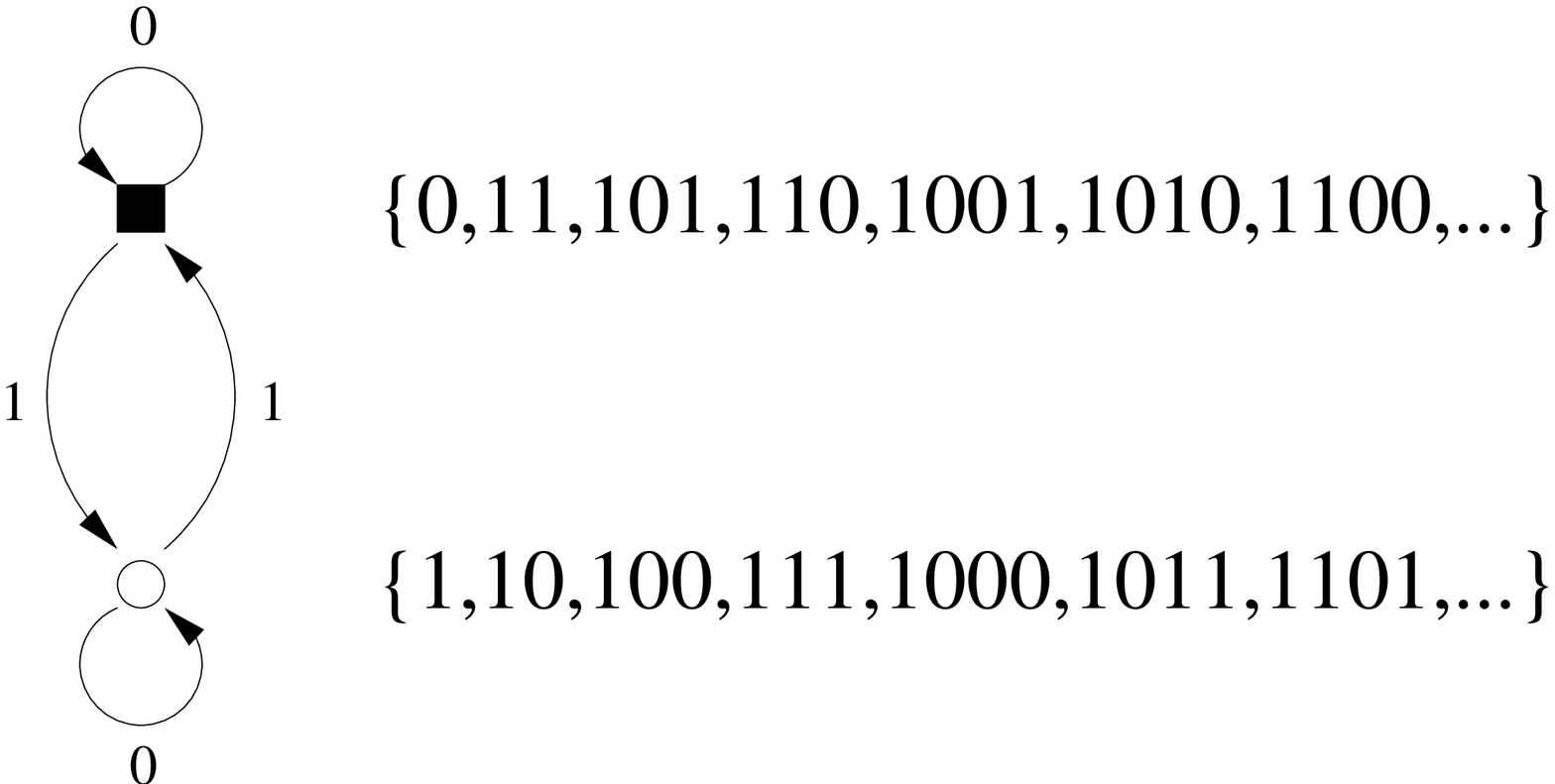}}
\end{example}
\begin{example}
Let $p=2$ and let $S$ be the
set of all nonnegative integers which have only $0$'s and $1$'s in their
base $4$ extension. 
We have
$$
S=\{0,1,4,5,16,17,20,21,\dots\}$$
and
$$
S=\{0,1,100,101,10000,10001,10100,10101,\dots\}
$$
in base 2.
The  automaton producing $S$ reversely  is as follows:\\[10pt]
\centerline{\includegraphics[scale=.4]{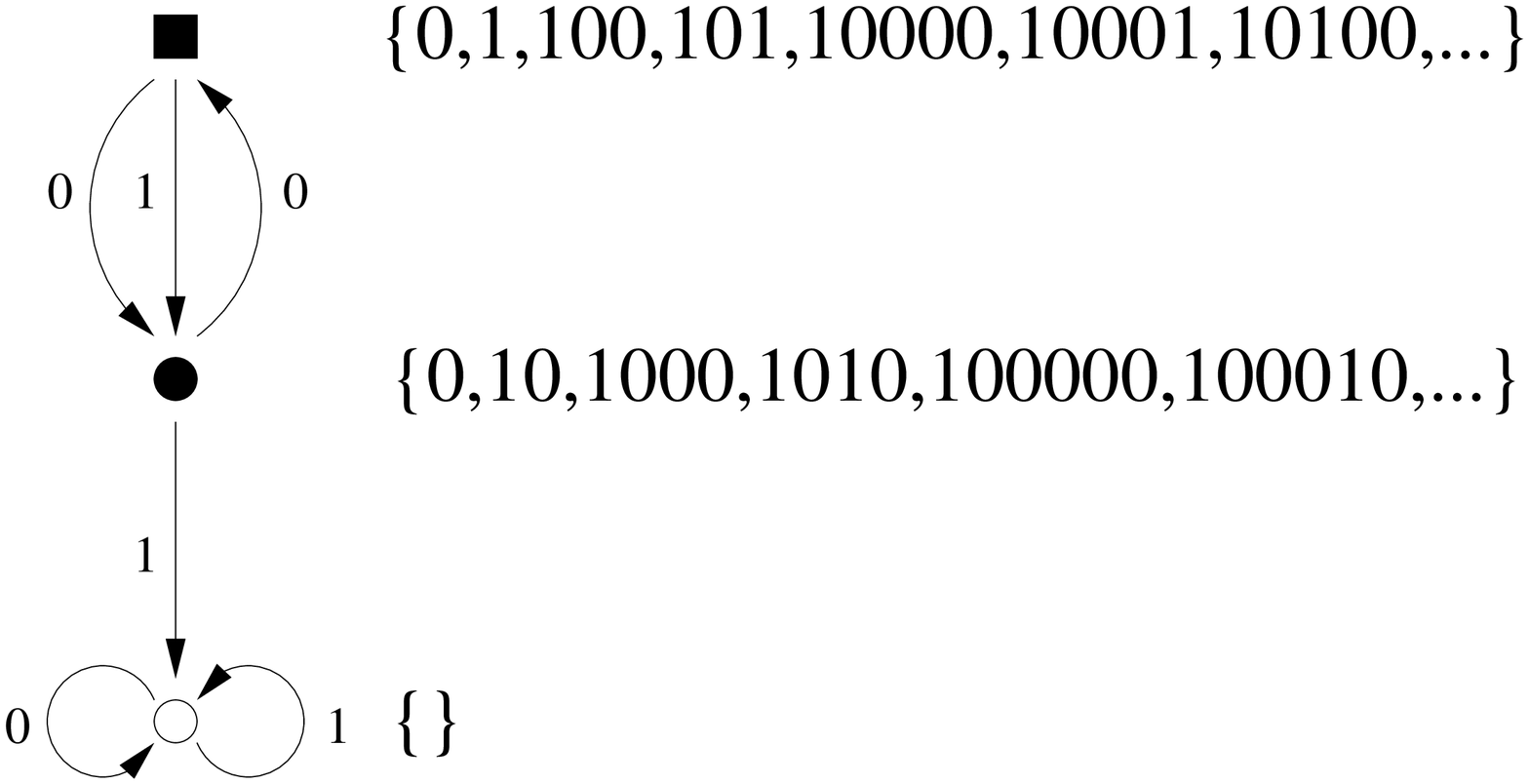}}
\end{example}
\begin{lemma}\label{lem7}
Suppose that $S\subseteq \N$ is a subset and $N$ is a positive integer. Then $S$ is $p$-automatic
if and only if $(L_j^N)^{-1}(S)$ is $p$-automatic for $j=0,1,2,\dots,N-1$.
\end{lemma}
\begin{proof}
The set $S$ is $p$-automatic if and only if
\begin{equation}\label{eqq1}
\{(L_j^{N})^{-1}(L_{i}^{p^e})^{-1}(S)\mid e\in \N, 0\leq i<p^e,0\leq j<N\}.
\end{equation}
is finite. On the other hand $(L_k^N)^{-1}(S)$ is $p$-automatic
for $k=0,1,\dots,N-1$ if and only if
\begin{equation}\label{eqq2}
\{(L_l^{p^e})^{-1}(L_k^N)^{-1}(S)\mid e\in \N, 0\leq l<p^e,0\leq k<N\}.
\end{equation}
is finite. If $p^ej+i=Nl+k$ with  $0\leq i,l<p^e$ and $0\leq j,k<N$ then we 
have
$$
L_i^{p^e}L_j^{N}=L^{Np^{e}}_{p^ej+i}=L^{Np^e}_{Nl+k}=L^{N}_kL^{p^e}_l.
$$
It follows that (\ref{eqq1}) and (\ref{eqq2}) are the same.
\end{proof}
The previous lemma shows that any arithmetic progression
is $p$-automatic. A major step toward the proof of Theorem~\ref{thm1a} is
the following result.
\begin{theorem}\label{thm2}
If $K$ is a field of characteristic $p>0$ and $a\in K^{\N}$ is
a recurrence sequence, then the ${\mathcal Z}(a)$
is $p$-automatic.
\end{theorem}
The next section is dedicated to the proof of Theorem~\ref{thm2}.
\section{Free Frobenius splitting}
\noindent Suppose that $K$ is a finitely generated field over $\F_p$.
For any subset $V$ of $K$ we define
$$
V^{\langle p \rangle}=\{f^p\mid f\in V\}.
$$
For $\alpha\in \R$ we denote the largest integer $\leq \alpha$
by $\lfloor \alpha\rfloor$.
\begin{lemma}\label{lem1}
Suppose that $V$ is an $n$-dimensional $\F_p$-subspace of $K$ containing $1$.
Then we have
$$
V^l\subseteq (V^{\lfloor l/p\rfloor})^{\langle p \rangle}V^{n(p-1)}.
$$
\end{lemma}
\begin{proof}
Let $e_1,\dots,e_n$ be a basis of $V$. Then $V^l$ is spanned
by monomials in $e_1,\dots,e_n$ of degree $l$. Let 
$$
f=e_1^{a_1}e_2^{a_2}\cdots e_n^{a_n}
$$
be such a monomial.
We can write $a_i=pb_i+c_i$ with $0\leq c_i\leq(p-1)$.
Then we have
$$
f=u^pv
$$
where 
$u=\prod_ie_i^{b_i}\in V^{\lfloor l/p\rfloor}$ and $v=\prod_i e_i^{c_i}\in V^{n(p-1)}$.
\end{proof}

Let $K^{\langle p \rangle}=\{f^p\mid f\in K\}$ be the subfield of $K$ of all $p$-th powers.
The field extension $K:K^{\langle p \rangle}$ has  finite degree, say $m$, and we can write
$$
K=K^{\langle p \rangle}h_1\oplus K^{\langle p \rangle}h_2\oplus \cdots\oplus K^{\langle p \rangle}h_m
$$
for certain $h_1,h_2,\dots,h_m\in K$.
Define $\pi_1,\pi_2,\dots,\pi_m:K\to K$ by
$$
f=\sum_{i=1}^m \pi_i(f)^ph_i
$$
for all $f\in K$. Note that $\pi_i(g^pf)=g\pi_i(f)$ for all $h,f\in K$.

\begin{proposition}\label{prop11}
Suppose that $V\subseteq K$ is a finite dimensional $\F_p$-subspace.
Then there exists a finite dimensional $\F_p$-subspace $W$
of $K$ containing $V$ such that
$$
\pi_i(VW)\subseteq W
$$
for all $i$.
\end{proposition}
\begin{proof}
Without loss of generality we way assume that $V$ contains
$1,h_1,\dots,h_m$ and generators of the field $K$ over $\F_p$.
Let $R$ be the ring generated by $V$. We have
$$
R\supseteq R^{\langle p \rangle}h_1\oplus R^{\langle p \rangle}h_2\oplus \cdots \oplus R^{\langle p \rangle}h_m.
$$
The module $M=R/(R^{\langle p \rangle}h_1\oplus R^{\langle p \rangle}h_2\oplus \cdots \oplus R^{\langle p \rangle}h_m)$
is a finitely generated torsion $R^{\langle p\rangle}$-module.
There exists a nonzero $g\in R$, such that $g^pM=0$. If we localize with respect
to $g$ we get
$$
R_g=R_{g^p}=(R_g)^{\langle p \rangle}h_1\oplus \cdots\oplus (R_g)^{\langle p
\rangle}h_m.
$$
Without loss of generality we may assume that $V$ contains $g^{-1}$.

The above discussion shows
that we only have to prove the proposition in the case where
the ring $R$ generated by $V$ satisfies
\begin{equation}\label{eqdirsum}
R=R^{\langle p \rangle}h_1\oplus \cdots \oplus R^{\langle p \rangle}h_m.
\end{equation}
Let $n:=\dim_{\F_p}V$. By Lemma~\ref{lem1} we have
\begin{equation}\label{eqpfloor}
V^l\subseteq (V^{\lfloor l/p\rfloor})^{\langle p \rangle}V^{n(p-1)}.
\end{equation}
From (\ref{eqdirsum}) follows that there exists a constant $C$ such that
\begin{equation}\label{eqpfloor2}
V^{n(p-1)}\subseteq (V^C)^{\langle p \rangle}h_1\oplus \cdots \oplus (V^C)^{\langle p \rangle}h_m.
\end{equation}
Combining (\ref{eqpfloor}) and (\ref{eqpfloor2}) gives
$$
V^l\subseteq (V^{C+\lfloor l/p\rfloor})^{\langle p \rangle}h_1\oplus \cdots
(V^{C+\lfloor l/p\rfloor})^{\langle p \rangle}h_m.
$$
If $l>Cp/(p-1)$ then
$$
C+\lfloor\textstyle \frac{l}{p}\rfloor\leq
C+\frac{l}{p}<\frac{l(p-1)}{p}+\frac{l}{p}=l
$$
which implies that
$$
V^l\subseteq (V^{l-1})^{\langle p \rangle}h_1\oplus \cdots \oplus (V^{l-1})^{\langle p \rangle}h_m.
$$
So we may take $W=V^{l-1}$.
\end{proof}
\begin{proof}[Proof of Theorem~\ref{thm2}]
If we change finitely many entries in the sequence $a$ then
${\mathcal Z}(a)$ will stay the same up to a finite set.
Without loss of generality we may assume that $a$ is basic by changing finitely many
entries in $a$.
In regard of Lemma~\ref{lem7} and Lemma~\ref{lem3.4} we may also assume that $a$ is
simple and nondegenerate.
 This means that we can write (after
enlarging $K$)
$$
a(n)=\sum_{i=1}^d\beta_i\alpha_i^n
$$
with $\alpha_1,\dots,\alpha_d,\beta_1,\dots,\beta_d\in K\setminus\{0\}$.
We may assume without loss of generality that $\alpha_1,\dots,\alpha_d$
are distinct.
Let $V\subseteq K$ be the $\F_p$-space spanned by 
all $\alpha_i^j$ with $0\leq j<p$ and all $\beta_i$.
Let $W$ be as in the Proposition~\ref{prop11}.
Define $a_i=(\alpha_i^0,\alpha_i^1,\dots)$.
Then we can write
$$
a=\beta_1a_1+\beta_2a_2+\cdots+\beta_ra_d.
$$
Consider the $\F_q$-vector space
$$
U=Wa_1+Wa_2+\cdots+Wa_d\subseteq W^{\N}
$$
(this is in fact a direct sum).
We claim that
$$
\pi_i(T^p_{j}U)\subseteq U
$$
for all $i\in \{1,2,\dots,m\}$ and all $j\in \{0,1,\dots,p-1\}$.
We have
$$
\pi_iT^p_j(Wa_k)=
\pi_iT^p_j(W(\alpha_k^0,\alpha_k^1,\dots))=
\pi_i(W(\alpha_k^j,\alpha_k^{j+p},\alpha_k^{j+2p},\dots))=
$$
$$
=\pi_i((W\alpha^j)(\alpha_k^0,\alpha_k^p,\dots))\subseteq 
\pi_i(WV)(\alpha_k^0,\alpha_k^1,\alpha_k^2,\dots)\subseteq
Wa_k.
$$
This proves that $\pi_i(T_j^pU)\subseteq U$.

Let ${\mathcal W}$ be the set of all
$$
{\mathcal Z}(b_1)\cap \cdots \cap {\mathcal Z}(b_r)
$$
with $r\in \N$ and $b_1,\dots,b_r\in U$. If $S\in {\mathcal W}$, then
$$
S={\mathcal Z}(c_1)\cap \cdots \cap{\mathcal Z}(c_t)
$$
for some $c_1,\dots,c_t\in U$.
Now we get
$$
(L_{i}^{p})^{-1}(S)=\bigcap_{j=1}^t (L_i^p)^{-1}({\mathcal Z}(c_j))
$$
and
$$
(L^p_i)^{-1}({\mathcal Z}(c_j))={\mathcal Z}(T^p_ic_j)=
{\mathcal Z}(\pi_1(T^p_ic_j))\cap \cdots \cap
{\mathcal Z}(\pi_m(T^p_ic_j))\in {\mathcal W}
$$
for all $j$ because $\pi_l(T^p_ic_j)\in U$ for all $l$ and $j$.
We obtain $(L_{i}^p)^{-1}(S)\in {\mathcal W}$ because
${\mathcal W}$ is closed under intersections.
This shows that ${\mathcal W}$ is closed under
the operations $S\mapsto (L^p_i)^{-1}(S)$ for $i=0,1,\dots,p-1$.
Because ${\mathcal Z}(a)\in {\mathcal W}$ we get
$${\mathcal V}_{{\mathcal Z}(a)}\subseteq {\mathcal W}
$$
Since ${\mathcal W}$ is finite, so is ${\mathcal V}_{{\mathcal Z}(a)}$.
We conclude that ${\mathcal Z}(a)$ is $p$-automatic.
\end{proof}
\section{Bounds for zero sets}
In this section we find explicit bounds for the zero set of recurrence sequences
in positive characteristic. The following sections do not depend on this section,
so this section may be skipped.
\begin{definition}
Suppose that $S\subseteq \N$ is a $p$-automatic set. We define
the {\it $p$-complexity of $S$} by
$$\comp_p(S):=|{\mathcal V}_S|,$$ 
where ${\mathcal V}_S$
is as in Corollary~\ref{corVS}. 
\end{definition}
The $p$-complexity  is a useful measure for complexity of a $p$-automatic
set as the following lemmas show.
\begin{lemma}\label{lemmm}
If $m\in S$ and $S$ is finite, then $m< p^{\comp_p(S)-2}$.
\end{lemma}
\begin{proof}
Without loss of generality we may assume that $m$ is the largest
element of $S$. If $m=0$ then $S=\{0\}$ and ${\mathcal
V}_S=\{\{0\},\emptyset\}$. In this case $\comp_p(S)=2$ and 
$$
m=0 <1=p^0=p^{\comp_p(S)-2}.
$$
Suppose that $m=[w_{r}\cdots w_0]_p$ with $w_0,\dots,w_r\in {\mathcal
A}:=\{0,1,\dots,p-1\}$ and
$w_r\neq 0$. The largest element
of $L_{w_0}^{-1}(S)$ is $[w_r\cdots w_1]_p$. Continuing in
this way we say that ${\mathcal V}_S$ contains a set whose
largest element is $[w_r\cdots w_s]_p$ for $s=0,1,\dots,r$.
The set ${\mathcal V}_S$ also contains $\{0\}$ and $\emptyset$.
So we have $\comp_p(S)=|{\mathcal V}_S|\geq r+3$.
On the other hand, $m<p^{r+1}$. We have
$$
m<p^{r+1}=p^{(r+3)-2}\leq p^{\comp_p(S)-2}.
$$
\end{proof}
\begin{lemma}\label{lemfinitebound}
Suppose that $S\subseteq \N$ is $p$-automatic and nonempty.
If $m\in S$ is the smallest element, then
 $m< p^{\comp_p(S)-2}$.
 \end{lemma}
 \begin{proof}
 The proof goes the same as for the previous lemma.
 \end{proof}
 \begin{lemma}\label{lembound}
 Suppose that $S_1,S_2\subseteq \N$ are both $p$-automatic.
 Then we have
 $$
 \comp_p((S_1\setminus S_2)\cup (S_2\setminus S_1))\leq
 \comp_p(S_1)\comp_p(S_2).
 $$
 \end{lemma}
 \begin{proof}
 Every element of ${\mathcal V}_{(S_1\setminus S_2)\cup (S_2\setminus S_1)}$ is of the form
$$
(U_1\setminus U_2)\cup (U_2\setminus U_1)
$$
with $U_1\in {\mathcal V}_{S_1}$ and $U_2\in {\mathcal V}_{S_2}$.
\end{proof}
Similar proofs show that 
$$\comp_p(S_1\cap S_2)\leq \comp_p(S_1)\comp_p(S_2),$$
$$\comp_p(S_1\cup S_2)\leq \comp_p(S_1)\comp_p(S_2),$$ 
and so forth.

\begin{proposition}\label{propbound}
Suppose that $\alpha_1(x),\dots,\alpha_d(x),\beta_1(x),\dots,\beta_d(x)\in 
\F_q[x]$ where $q$ is a power of the prime $p$
and define $a\in \F_q[x]^{\N}$ by
$$
a(n)=\beta_1(x)\alpha_1(x)^n+\cdots+\beta_d(x)\alpha_d(x)^n.
$$
Suppose that $k\geq 2$ such that
$\deg(\alpha_i(x))\leq k$ and $\deg(\beta_i(x))\leq (p-1)k$
for all $i$. Then we have 
$$\comp_p({\mathcal Z}(a))\leq 
q^{d^2k^4p^4}.
$$
\end{proposition}
\begin{proof}
Let $V$ be the space all polynomials of degree $\leq (p-1)k$.
We choose $W$ as in Proposition~\ref{prop11}. We follow the proof
of Proposition~\ref{prop11} to find $W$ explicitly. Let $m:=p$ and define
$h_1,\dots,h_p$ by $h_i=x^{i-1}$. We get a decomposition
$$
\F_q[x]=\F_q[x]h_1\oplus \cdots\oplus \F_q[x]h_p.
$$
Let $n=\dim_{\F_q}V=k(p-1)+1$.
We need to choose $C$ such that
\begin{equation}\label{eqsubset}
V^{n(p-1)}\subseteq 
(V^C)^{\langle p\rangle}h_1\oplus \cdots\oplus (V^C)^{\langle p\rangle} h_p.
\end{equation}
Now $V^{n(p-1)}$ is the set of
all polynomials of degree $\leq n(p-1)^2k$ and
$$
(V^C)^{\langle p\rangle}h_1\oplus \cdots\oplus (V^C)^{\langle p\rangle} h_p
$$
is the set of all polynomials of degree at most $Cpk(p-1)+(p-1)$.
It suffices that
$$
n(p-1)^2k\leq Cpk(p-1)+(p-1)
$$
which is equivalent to
$$
n(p-1)k\leq Cpk+1
$$
and (because $k\geq 2$) to
$$n\leq \frac{Cp}{p-1}.$$
We  take $C=\lceil n(p-1)/p\rceil$, so that
\begin{multline}\label{eqCp}
\frac{Cp}{p-1}\leq
\frac{n(p-1)+(p-1)}{p}=\frac{(n+1)(p-1)}{p}=\\
=\frac{((p-1)k+2)(p-1)}{p}<
\frac{k(p+1)(p-1)}{p}<kp.
\end{multline}
In the proof of Proposition~\ref{prop11} we must take $l$ such that
$l>Cp/(p-1)$. So we may take $l=kp$ by (\ref{eqCp}).
Define $W=V^{l-1}=V^{kp-1}$ which
is the set of polynomials of degree $\leq (kp-1)k(p-1)$.
We have $\dim W=(kp-1)k(p-1)+1\leq k^2p^2$.

For $U$ as in the proof of Theorem~\ref{thm2} we have
$$
u:=\dim U \leq d\cdot \dim W=dk^2p^2.
$$
Let ${\mathcal W}$ be the set of all 
\begin{equation}\label{eqint}
{\mathcal Z}(b_1)\cap \cdots\cap {\mathcal Z}(b_r)
\end{equation}
where $r\in \N$ and $b_1,\dots,b_r\in U$.
Note that (\ref{eqint}) only depends on the $\F_q$-vector space
spanned by $b_1,\dots,b_r$. Therefore $|{\mathcal W}|$ is
bounded by the number of subspaces of $U$. 
Every subspace of $U$ can be generated by a $u\times u$ matrix
with entries in $\F_q$. So
a (rough) upper bound for the number of subspaces of $U$ is
$$
q^{u^2}\leq q^{d^2k^4p^4}.
$$
We conclude that
$$
|{\mathcal V}_{{\mathcal Z}(a)}|\leq |{\mathcal W}|\leq q^{d^2k^4p^4}.
$$
\end{proof}
\begin{corollary}
In the setup of Proposition~\ref{propbound}, we have
then
$$
\min({\mathcal Z}(a))< p^{\displaystyle q^{d^2k^4p^4}}.
$$
\end{corollary}
\begin{proof}
This follows from Proposition~\ref{propbound} and Lemma~\ref{lemfinitebound}.
\end{proof}
\begin{corollary}
In the setup of Proposition~\ref{propbound}, if ${\mathcal Z}(a)$ is finite,
then
$$
\max({\mathcal Z}(a))< p^{\displaystyle q^{d^2k^4p^4}}.
$$
\end{corollary}
\begin{proof}
This follows from Proposition~\ref{propbound} and
Lemma~\ref{lemmm}.
\end{proof}
\begin{proposition}\label{prop6.7}
Suppose that $K$ is a finitely
generated field over $\F_p$ and
 $a\in K^{\N}$ is a $K$-recurrence sequence
 which is simple and nondegenerate. (Assume that the recurrence
 relation is explicitly known.)
 One can compute
an explicit bound $N(a)$ such that
$$
\comp_p({\mathcal Z}(a))\leq N(a).
$$
\end{proposition}
\begin{proof}
 By possibly
enlarging the field $K$ we can explicitly write
$$
a(n)=\beta_1\alpha_1^n+\cdots+\beta_d\alpha_d^n
$$
where $\alpha_1,\dots,\alpha_d,\beta_1,\dots,\beta_d\in K$.
As in the proof of Proposition~\ref{propbound}, we can follow the proofs
of Proposition~\ref{prop11} and Theorem~\ref{thm2} to find
an explicit bound for $N(a)$.
\end{proof}
\begin{proof}[Proof of Theorem~\ref{theoremeffective}]
We can reduce to the case where $a$ is simple and
nondegenerate by Lemma~\ref{lem3.4}.
It is possible to explicitly enumerate all
$p$-normal sets $S_1,S_2,\dots$. 
We can verify whether ${\mathcal Z}(a)=S_i$ as follows.
Proposition~\ref{prop6.7}
gives an upper bound for $\comp_p({\mathcal Z}(a))$. One can
explicitly construct an automaton that produces $S_i$ reversely.
This gives an upper bound $\comp_p(S_i)$.
Let
$$
U_i=({\mathcal Z}(a)\setminus S_i)\cup (S_i\setminus{\mathcal Z}(a)).
$$
Then $\comp_p(U_i)\leq \comp_p({\mathcal Z}(a))\comp_p(S_i)$ by 
Lemma~\ref{lembound}, so
we have an explicit upper bound for $\comp_p(U_i)$, say
$\comp_p(U_i)\leq N$.
If $U_i$ is nonempty then the smallest element of $U_i$ is at most
$p^{N-2}$ (see Lemma~\ref{lemfinitebound}). So we have that $S_i={\mathcal Z}(a)$ if and only if
$$
S_i\cap \{0,1,2,\dots,p^{N-2}\}={\mathcal Z}(a)\cap \{0,1,2,\dots,p^{N-2}\}
$$
and this can be verified in a finite amount of time.
\end{proof}


\section{Automata producing zero sets}
For any simple nondegenerate linear recurrence sequence $a\in K^\N$ where
$K$ is a field of characteristic $p>0$ we constructed
a  $p$-automaton that produces ${\mathcal Z}(a)$ reversely. In this
section we will study how such an automaton can look like.
As it turns out, these automata have a very special form.

\begin{definition}
Suppose that $S\subseteq \N$ is of the form ${\mathcal Z}(a)$
for some simple nondegenerate recurrence sequence $a\in K^{\N}$ where $K$ is a field
of characteristic $p>0$.
We define the {\it level\/} $\ell(S)$ of $S$ as the smallest nonnegative
integer $d$ for which we can write
$$
S={\mathcal Z}(b_1)\cap \cdots \cap {\mathcal Z}(b_r)
$$
where $b_1,\dots,b_r\in K^{\N}$ are linear recurrence sequences of
order $\leq d$.
\end{definition}
\noindent For example $\ell(\N)=0$ because
the zero sequence has minimum polynomial $1$. 
Also  $\ell(\emptyset)=1$ because the sequence that is constant $1$
has minimum polynomial $S-1$. If $S\neq \emptyset,\N$ then $\ell(S)\geq 2$
because any linear recurrence sequence of order 1 is  constant.
\begin{example}
Let $a\in \F_2(x)^{\N}$ be defined by
$$
a(n)=(x+1)^n-x^n-1
$$
(See Example~\ref{ex2}). The
 automaton producing ${\mathcal Z}(a)=\{1,2,4,8,\dots\}$ reversely is:\\[10pt]
\centerline{\includegraphics[scale=.4]{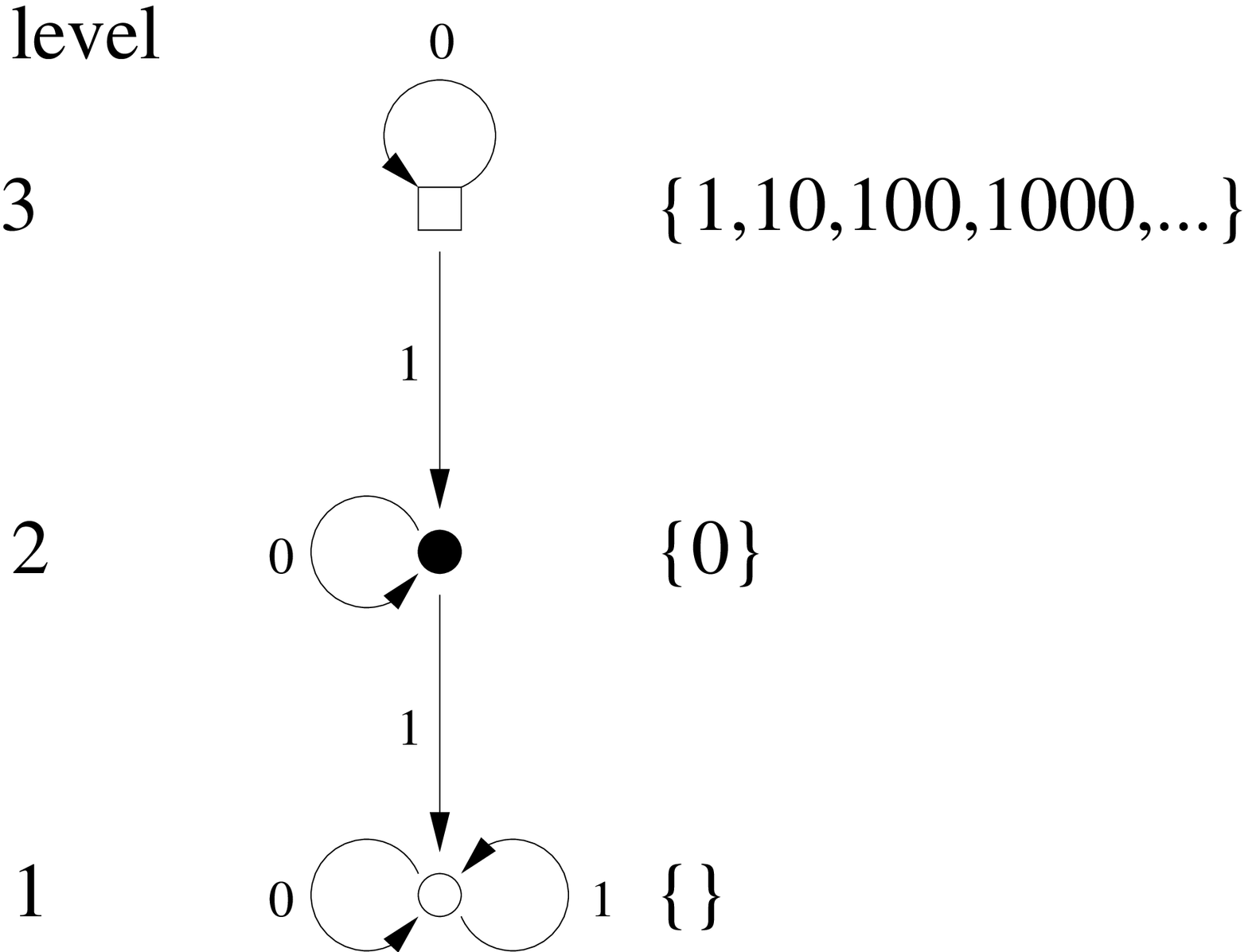}}
(elements of $\N$ are written in base 2).
\end{example}
\begin{example}
Let $a\in \F_3(x)^{\N}$ be defined by
$$
a(n)=(x+1)^n-x^n-1
$$
(See Example~\ref{ex2}). The
 automaton producing 
${\mathcal Z}(a)=\{1,3,9,27,\dots\}$ reversely is:\\[10pt]
\centerline{\includegraphics[scale=.4]{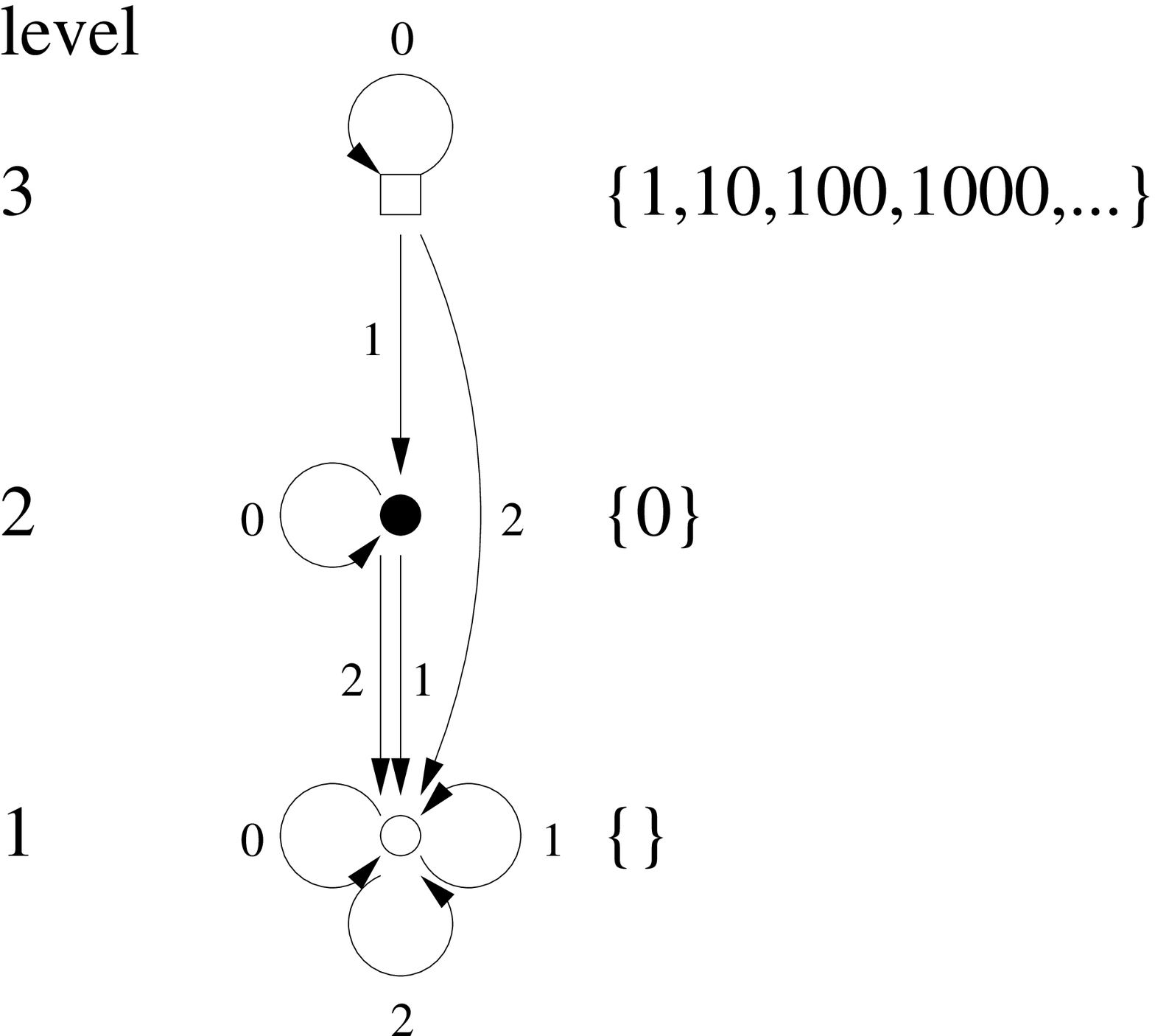}}
(elements in $\N$ are written in base 3).
\end{example}
\begin{example}
Let $a\in \F_2(x,y,z)^{\N}$ be defined by
$$
a(n)=(x+y+z)^n-(x+y)^n-(x+z)^n-(y+z)^n+x^n+y^n+z^n
$$
(See Proposition~\ref{prop1}). The  automaton producing 
${\mathcal Z}(a)=\{2^i+2^j\mid i,j\in \N\}\cup \{1\}$ reversely is:\\[10pt]
\centerline{\includegraphics[scale=.4]{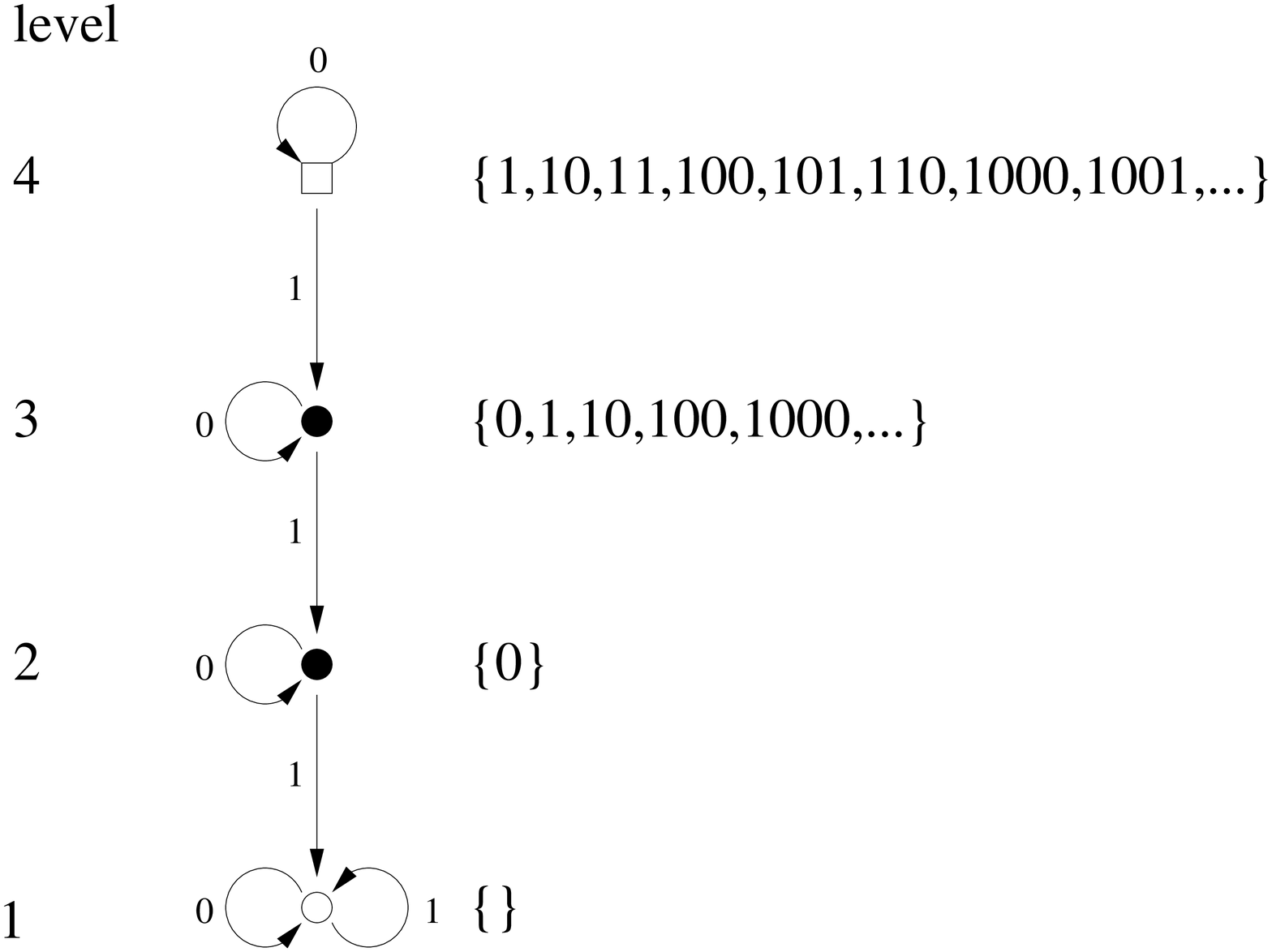}}
(in base $2$).
\end{example}
\begin{example}
Let $a\in \F_4(x)^{\N}$ be defined by
$$
a(n)=x^n+(x+1)^n+(x+\alpha)^n+(x+1+\alpha)^n
$$
as in Example~\ref{ex6}. The  automaton producing ${\mathcal Z}(a)$ reversely is:\\[10pt]
\centerline{\includegraphics[scale=.4]{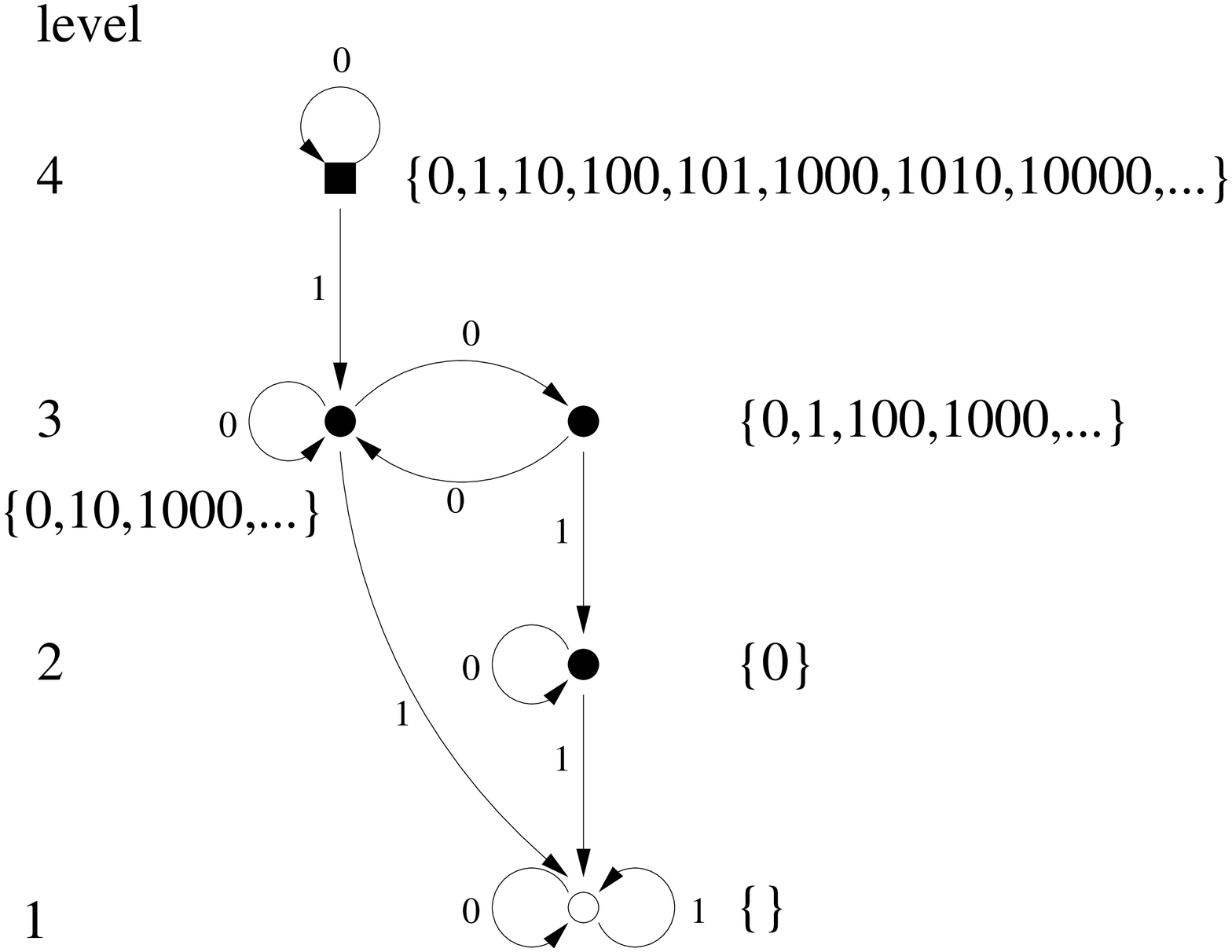}}
(in base 2).
\end{example}

\begin{proposition}\label{prop5.3}
Suppose that $a$ is a simple nondegenerate nonzero recurrence sequence of order $d$ in
a field $K$ of characteristic $p>0$. Consider the automaton that produces
$S={\mathcal Z}(a)$ reversely and its graph.
\begin{enumerate}
\letters
\item If $Q,R\in {\mathcal V}_S$ and there is a path 
from $Q$ to $R$ then $\ell(Q)\geq \ell(R)$.
\item If $Q,R\in {\mathcal V}_S$ with $\ell(Q)\geq 2$ 
and there are two distinct paths
 from $Q$ to $R$ of the same length, then $\ell(Q)>\ell(R)$.
\item $\N\not\in {\mathcal V}_S$.
\item  $\emptyset \in {\mathcal V}_S$.
\end{enumerate}
\end{proposition}
\begin{proof}
(a) Suppose that $R=(L_i^p)^{-1}(Q)$ and that
$$
Q=\bigcap_{j=1}^r{\mathcal Z}(b_j)
$$
where $b_j$ is a simple nondegenerate 
recurrence sequence of order $\leq \ell(Q)$ for all $j$.
Then 
$$
R=\bigcap_{j=1}^r{\mathcal Z}(T_i^pb_j)
$$
and $T_i^pb_j$ has order $\leq \ell(Q)$ for all $j$ by Lemma~\ref{lem2.1}. It
follows that $\ell(R)\leq \ell(Q)$.

(b) Suppose that there are 2 paths from $Q$ to $R$ of length $r$.
We can write
$$
Q={\mathcal Z}(a_1)\cap \cdots \cap {\mathcal Z}(a_s)
$$
where $a_i$ is a simple nondegenerate recurrence sequence of order at most $d:=\ell(Q)$.
Since there are two paths of length $r$
from $Q$ to $R$, we have  
$$(L_j^{p^r})^{-1}(Q)=R=(L_k^{p^r})^{-1}(Q)$$
for some $j,k$ with $0\leq j<k<p^r$.
This implies
$$
R=\bigcap_{i=1}^r{\mathcal Z}(T^{p^r}_ja_i)
$$
and
$$
R=\bigcap_{i=1}^r {\mathcal Z}(T^{p^r}_ka_i)
$$
By Lemma~\ref{lem5} below we can write
$$
{\mathcal Z}(T^{p^r}_ja_i)\cap {\mathcal Z}(T^{p^r}_ka_i)=
{\mathcal Z}(b_i)\cap {\mathcal Z}(c_i)
$$
for certain recurrence sequences $b_i$ and $c_i$ of order $\leq d-1$.
We see that
$$
R=\bigcap_{i=1}^r\big({\mathcal Z}(T^{p^r}_ja_i)\cap {\mathcal
Z}(T^{p^r}_k a_i)\big)=\bigcap_{i=1}^r\big(
{\mathcal Z}(b_i)\cap {\mathcal Z}(c_i)\big).
$$
It follows that $\ell(R)\leq d-1<\ell(Q)$.

(c) For every $Q\in {\mathcal V}_S$ there exists an $r$ and $j$
such that $Q={\mathcal Z}(T^{p^r}_ja)$. Since $a$ is simple and nondegenerate,
$T^{p^r}_ja$  cannot be the 0 sequence by Lemma~\ref{lem2.6}.

(d) Let $Q\in {\mathcal V}_S$ with $\ell(Q)$ minimal. There are $p^r$
paths starting at $Q$ of length $r$. Choose $r$ such that $p^r>|{\mathcal V}_S|$.
By the pigeonhole principle, there are two paths of length $r$
which have the same endpoint, say $R$. If $\ell(Q)\geq 2$
then $\ell(Q)>\ell(R)$ by part (b) and we have a contradiction. Therefore,
$\ell(Q)\leq 1$, so  $Q=\emptyset$ or $Q=\N$. But $Q\neq \N$ by part (c).
\end{proof}
\begin{lemma}\label{lem5}
Suppose that $a\in K^{\N}$ is a simple and nondegenerate sequence of order 
$d\geq 2$ where $K$
is a field of characteristic $p>0$. 
If  $j\neq k$ then there exist simple and nondegenerate sequences $b,c\in K^{\N}$
of order $\leq d-1$ such that
$$
{\mathcal Z}(T^{p^r}_ja)\cap {\mathcal Z}(T^{p^r}_ka)={\mathcal Z}(b)\cap {\mathcal Z}(c).
$$
\end{lemma}
\begin{proof}
We can write 
$$
a(n)=\sum_{i=1}^d\beta_i\alpha_i^n.$$
We have
$$
(T^{p^r}_ja)(n)=\sum_{i=1}^d(\beta_i\alpha_i^j)(\alpha_i^{p^r})^n
$$
(and a similar formula for $T^{p^r}_ka$).
Define
$$
b=\alpha_1^kT^{p^r}_ja-\alpha_1^jT^{p^r}_ka
$$
and
$$
c=\alpha_2^kT^{p^r}_j a-T^{p^r}_j\alpha_2^jT^{p^r}_ka.
$$
By construction, the coefficient of $(\alpha_1^{p^r})^n$
in $b(n)$ vanishes. Similarly, the coefficient
of $(\alpha_2^{p^r})^n$ in $c(n)$ vanishes. This shows
that $b$ and $c$ have order $\leq d-1$.
Since $\alpha_1/\alpha_2$ is not a root of unity, we have
$$
\det\begin{pmatrix}
\alpha_1^k & -\alpha_1^j\\
\alpha_2^k & -\alpha_2^j\end{pmatrix}\neq 0.
$$
Therefore, the $K$-span of $b$ and $c$ is the same as the span of $T^{p^r}_ka$
and $T^{p^r}_ja$, hence
$$
{\mathcal Z}(T^{p^r}_ja)\cap {\mathcal Z}(T^{p^r}_ka)={\mathcal Z}(b)\cap {\mathcal Z}(c).
$$
\end{proof}

\begin{definition}
Let ${\mathcal A}=\{0,1,\dots,p-1\}$.
Suppose that $u_0,u_1,\dots,u_d,w_1,\dots,w_{d}\in{\mathcal A}^\star$. 
We define
$$U_p(u_0,u_1,\dots,u_d;w_1,\dots,w_d)=
\{[u_dw_d^{k_d}u_{d-1}w_{d-1}^{k_{d-1}}\cdots u_1w_1^{k_1}u_0]_p\mid
k_1,k_2,\dots,k_d\in \N\}\subseteq\N.
$$
\end{definition}
\begin{proposition}\label{prop6.9}
Let $a\in K^\N$ be a nonzero, simple and nondegenerate sequence of order $d$, where $K$ is a field
of characteristic $p>0$. Then $S:={\mathcal Z}(a)$ is a finite
union of sets of the form
$U_p(u_0,\dots,u_m;w_1,\dots,w_m)$ with $m\leq d-2$.
\end{proposition}
\begin{proof}
 If $w=i_ri_{r-1}\cdots i_0$
then the path given by the word $w$ starting at some vertex $Q$
is the path in the graph that visits the vertices
$$
Q,i_0\cdot Q,i_1i_0\cdot Q,\dots,i_ri_{r-1}\cdots i_0\cdot Q.
$$

A vertex $Q$ in ${\mathcal V}_S$ is
called a {\it loop vertex\/} if there is a nontrivial path from $Q$ to itself.
In other words, $Q$ is a loop vertex if
and only if $w\cdot Q=Q$ for some nontrivial word $w$. 
Let us choose $w$ nontrivial of minimal length such that $w\cdot Q=Q$.
Then the path from $Q$ to $Q$ via $w$ does not intersect itself
(except at the beginning and the end). We claim that every path
from $Q$ to $Q$ is given by a power of $w$. Suppose that $u\cdot Q=Q$.
We can write $u=yw^s$ where $y$ is a path that does not have $w$
as a prefix, i.e., $y$ is not of the form $zw$ for some word $z$. 
If $y$ is the trivial path then were are done, so assume to the
contrary that $y$ is not trivial.
Let $e$ and $f$ be the lengths of $w$ and $y$ respectively.
Then $y^e$ and $w^f$ are paths from $Q$ to $Q$  of length $ef$.
Since tautologically $\ell(Q)=\ell(Q)$, we have $y^e=w^f$ by
Proposition~\ref{prop5.3}(b). Since $e\leq f$, $w$ has to be
a prefix of $y$. Contradiction. We conclude that every loop from $Q$
to $Q$ is given by a power of $w$.

Suppose that $n\in S$. We can write $n=[w]_p$ for some word $w$.
Consider the path $\gamma$ from $S$ to $Q:=w\cdot S$ given by $w$.
Define $S_0,S_1,S_2,\dots$ as follows. 
First we define $S_0=S$. For $j>0$ we define $S_j$ as
the first loop vertex in the path $\gamma$ from which there exists no
path to $S_{j-1}$ if such a vertex exists. Suppose that
we can define $S_1,S_2,\dots,S_m$ in this way.
 We define
$S_{m+1}=Q$. 
We can write
$$
w=
u_mw_m^{l_m}u_{m-1}w_{m-1}^{l_{m-1}}u_{m-2}\cdots u_1w_1^{l_1}u_0
$$
where $u_j$ defines a path from $S_{j}$ to $S_{j+1}$ 
without self-intersection,  $w_j$ defines the unique nontrivial loop at
vertex $S_{j}$ without self intersection, and $l_j\in \N$ for all $j$.

For $j$ with $1\leq j<m$, for every $i,k$ we have that
$$
w_{j+1}^lu_{j}w_{j}^k
$$
defines a path from $S_j$ to $S_{j+1}$. 
Let $s_j$ and $s_{j+1}$ be the lengths of $w_{j}$ and $w_{j+1}$ respectively.
Then $w_{j+1}^{s_j}u_j$ and $u_jw_{j}^{s_{j+1}}$ define paths
from $S_j$ to $S_{j+1}$ of equal length. Since $w_{j+1}^{s_j}u_j$
visits vertex $S_j$ only once at the beginning, and
$u_jw_j^{s_{j+1}}$ visits $S_j$ exactly $s_{j+1}+1$ times, we have
that
$w_{j+1}^{s_j}u_j\neq u_jw_j^{s_{j+1}}$.  From Proposition~\ref{prop5.3}(b)
follows that $\ell(S_{j+1})>\ell(S_j)$.

We deduce that  $S$ contains $U_p(u_0,u_1,\dots,u_m;w_1,\dots,w_m)$.
Note that
$$
d\geq \ell(S)=\ell(S_0)> \ell(S_1)>\ell(S_2)>\cdots>\ell(S_m)\geq \ell(S_{m+1})=
\ell(Q)\geq 2.
$$
It follows that $m\leq d-2$.

Since there are only finitely many paths without self-intersection
and only finitely many loops without self-intersection, it follows
that $S$ is a finite union of sets of the form
$U_p(u_0,u_1,\dots,u_m;w_1,\dots,w_m)$ with $m\leq d-2$.
\end{proof}
\section{Proof of the main result}
In this section will prove Theorem~\ref{thm1c}. This also completes
the proof of Theorem~\ref{thm1a}.
\begin{lemma}\label{lem5.8}
Suppose that $a\in K^\N$ is a simple nondegenerate sequence 
in a field $K$ of characteristic $p>0$. Suppose that $q$ is a 
power of $p$ and that $r,s\in \Q$.
If
$$
\{r+sq^n\mid n\geq m\}\subseteq {\mathcal Z}(a)
$$
for some constant $m\in \N$ then
$$
\{r+sq^n\mid n\in \N\}\cap \N\subseteq {\mathcal Z}(a).
$$
\end{lemma}
\begin{proof}
We can write
$$
a(n)=\sum_{i=1}^d\beta_i\alpha_i^n
$$
where $\alpha_i/\alpha_j$ is not a root of unity for all $i\neq j$.
Choose $N\in \N$ such that $Nr,Ns\in \Z$. Choose $\lambda_i$  such that
$\lambda_i^N=\alpha_i$. Then
 $\lambda_i/\lambda_j$ is not a root of unity for $i\neq j$
and
$$
b(n)=\sum_{i=1}^d\beta_i\lambda_i^{n}
$$
is a simple nondegenerate sequence with $b(Nn)=a(n)$
and 
$$b(Nr+(Ns)q^n)=a(r+sq^n).$$
The sequence $c(n)$ defined by 
$$c(n)=b(Nr+(Ns)n)=a(r+sn)$$ 
is also a
simple nondegenerate sequence.
So we can reduce the lemma to the case that $r=0$ and $s=1$.

Let us assume that 
$$
a(n)=\sum_{i=1}^d\beta_i\alpha_i^n.
$$
is a simple nondegenerate sequence with $a(q^n)=0$ for $n\geq m$. 
We will show that $a(q^n)=0$ for all $n\in\N$.
Consider the element
$$
x=\sum_{i=1}^d\beta_i\otimes \alpha_i\in K\otimes_{\F_q}K
$$
Suppose that $x\neq 0$.
Let us write
$$
x=\sum_{i=1}^e\delta_i\otimes\gamma_i
$$
where $e$ is minimal. This implies that $\gamma_1,\dots,\gamma_e$ are
linearly independent over $\F_q$. Similarly,
 $\delta_1,\dots,\delta_e$ are linearly independent over $\F_q$.
 Let $\phi:K\to K$ be the Frobenius homomorphism
 defined by $\phi(\epsilon)=\epsilon^q$.
 The homomorphism $\phi$ leaves the field $\F_q$ invariant.
 Define $\psi:K\otimes_{\F_q}K\to K$ by $\psi(\sum_i \lambda_i\otimes\nu_i)=
 \sum_i\lambda_i\nu_i$.
 We have
 $$
 \psi\circ (\id\otimes \phi^n)\big(\sum_{i=1}^d\beta_i\otimes \alpha_i\big)=
 \psi\big(\sum_{i=1}^d\beta_i\otimes \alpha_i^{q^n}\big)= \sum_{i=1}^d \beta_i\alpha_i^{q^n}=a(q^n).
$$
By assumption this is equal to $0$ for $n\geq m$.
On the other hand, this is equal to
$$
\psi\circ (\id\otimes \phi^n)\big(\sum_{i=1}^e\delta_i\otimes \gamma_i\big)=
 \psi\big(\sum_{i=1}^e\delta_i\otimes \gamma_i^{q^n}\big)=\sum_{i=1}^e \delta_i\otimes
 \gamma_i^{q^n}=
 \sum_{i=1}^e \delta_i\gamma_i^{q^n}
$$ 
Define $c_1,\dots,c_e\in K^\N$ by
$$
c_j(n)=\gamma_j^{q^n}.
$$
We know that $E^mc_1,\dots,E^mc_e$ are linearly dependent over $K$ because
$$
\Big(\sum_{j=1}^e \delta_jE^mc_j\Big)(n)=\sum_{j=1}^e\delta_j\gamma_j^{q^{m+n}}=a(q^{n+m})=0
$$
for all $n\in \N$.
Choose $j$ maximal such that
$E^nc_1,\dots,E^{n}c_{j-1}$ are linearly independent for all $n\in \N$.
Then $E^nc_1,\dots,E^nc_{j}$ are linearly dependent
for $n$ large enough.
There are unique $\epsilon_1,\dots,\epsilon_{j-1}\in K$
such that
\begin{equation}\label{eq20}
E^nc_j=\sum_{i=1}^{j-1}\epsilon_iE^nc_i.
\end{equation}
Taking the $q$-th power of (\ref{eq20}) yields:
\begin{equation}\label{eq21}
E^{n+1}c_j=(E^nc_j)^q=\sum_{i=1}^{j-1}\epsilon_i^q(E^nc_i)^q=
\sum_{i=1}^{j-1}\epsilon_i^qE^{n+1}c_i.
\end{equation}
Applying $E$ to (\ref{eq20}) yields:
\begin{equation}\label{eq22}
E^{n+1}c_j=\sum_{i=1}^{j-1}\epsilon_iE^{n+1}c_i.
\end{equation}
Subtracting (\ref{eq22}) from (\ref{eq21}) gives:
$$
0=\sum_{i=1}^{j-1}(\epsilon_i^q-\epsilon_i)E^{n+1}c_i
$$

Since $E^{n+1}c_1,E^{n+1}c_2,\dots,E^{n+1}c_{j-1}$ are linearly independent 
over $K$, we conclude that 
$$
\epsilon_i^q=\epsilon_i
$$
for $i=1,2,\dots,j-1$. This means that $\epsilon_1,\dots,\epsilon_{j-1}\in
\F_q$. Therefore 
$$E^nc_1=c_1^{q^n},\dots,E^nc_j=c_j^{q^n}$$ are 
linearly dependent over
$\F_q$. Taking the $q^n$-th root shows us that $c_1,\dots,c_j$ are
linearly dependent over $\F_q$. 
But then $\gamma_1,\dots,\gamma_j$ are linearly dependent over $\F_q$.
Contradiction! We conclude that $x=0$.
It follows that $a(q^n)=0$ for all $n$.
\end{proof}

\begin{proof}[Proof of Theorem~\ref{thm1c}]
By Proposition~\ref{prop6.9},
 ${\mathcal Z}(a)$ is
a finite union of sets of the form $U_p(u_0,u_1,\dots,u_m;w_1,\dots,w_m)$
with $m\leq d-2$. Let $r_i$ be the length of $w_i$ and
let $r$ be the least common multiple of $r_1,\dots,r_m$.
We can write $U_p(u_0,\dots,u_m;w_1,\dots,w_m)$ as a finite union
of sets of the form
$$U_p(u_0',\dots,u_m';w_1^{k_1},\dots,w_m^{k_m})$$
where $k_i=r/r_i$. This shows that ${\mathcal Z}(a)$ is 
 a finite union of sets of the form
$$U_p(u_0,u_1,\dots,u_m;w_1,\dots,w_m)$$ 
where $w_1,\dots,w_m$ all have the same
length $r$.
 Set $q=p^r$
and let $t_i$ be the length of $u_i$ for all $i$.
From
$$
[u_mw_m^{k_m}u_{m-1}w_{m-1}^{k_{m-1}}\cdots u_1w_1^{k_1}u_0]_p=
\sum_{i=0}^m [u_i]_pp^{t_0+t_1+\cdots+t_{i-1}}q^{k_1+k_2+\cdots+k_{i}}
+$$
$$+\sum_{i=1}^m
[w_i]_pp^{t_0+t_1+\cdots+t_{i-1}}q^{k_1+\cdots+k_{i-1}}\big(\frac{q^{k_i}-1}{q-1}\big)
$$
 follows that the set
$U_q(u_0,\dots,u_m;w_1,\dots,w_m)$
has the form
\begin{equation}\label{eqa}
\{c_0+c_1q^{l_1}+c_2q^{l_2}+\cdots+c_rq^{l_m}\mid 0\leq l_1\leq l_2\leq
\cdots\leq l_m\}
\end{equation}
with $c_0,c_1,\dots,c_m\in \Q$ such that $c_0+\cdots+c_m\in \Z$
and $(q-1)c_i\in \Z$ for all $i$.
We will show  that ${\mathcal Z}(a)$ contains
$$
\{c_0+c_1q^{l_1}+c_2q^{l_2}+\cdots+c_mq^{l_m}\mid l_1,\dots,l_m\in \N\}\cap \N.
$$
Suppose that $l_1,\dots,l_m\in \N$ such that
$$
c_0+c_1q^{l_1}+c_2q^{l_2}+\cdots+c_mq^{l_m}\in \N.
$$
We would like to show that
$$
c_0+c_1q^{l_1}+c_2q^{l_2}+\cdots+c_mq^{l_m}\in {\mathcal Z}(a).
$$
Since $c_m>0$, it suffices to show, by Lemma~\ref{lem5.8}, that
$$
c_0+c_1q^{l_1}+c_2q^{l_2}+\cdots+c_{m-1}q^{l_{m-1}}+c_mq^{l_m+N}\in {\mathcal Z}(a)
$$
for $N\gg 0$. We prove this by induction on 
$$
D:=|\{i\mid 1\leq i\leq m-1, l_i>l_{i+1}\}|.
$$
The case $D=0$ follows from (\ref{eqa}). If $D>0$, then there
exists an $i$ such that $l_i>l_{i+1}$. For $N$ sufficiently large
we have
$$
c_{i+1}q^{l_{i+1}}+c_{i+2}q^{l_{i+2}}+\cdots+c_{m-1}q^{l_{m-1}}+c_{m}q^{l_m+N}>0.
$$
For $M\geq l_{i}-l_{i+1}$ we get
$$
(c_1q^{l_1}+\cdots+c_{i}q^{l_i})+(c_{i+1}q^{l_{i+1}+M}+\cdots+
c_{m-1}q^{l_{m-1}+M}+c_{m}q^{l_m+M+N})\in {\mathcal Z}(a)
$$
by induction. From Lemma~\ref{lem5.8} follows that
$$
c_0+c_1q^{l_1}+c_2q^{l_2}+\cdots+c_{m-1}q^{l_{m-1}}+c_mq^{l_m+N}\in {\mathcal Z}(a)
$$
and we are done.
\end{proof}
We remark that Theorem~\ref{thm1c} is related to the following theorem
of Masser (see~\cite{Masser2}, \cite[\S 28]{Schmidt} and \cite[page 707]{KS}).
\begin{theorem}
Let $K$ be a field of characteristic $p> 0$ with algebraic
closure $\overline{K}\supseteq K$, $d\geq 1$ and let $\alpha_1,\dots,\alpha_d\in
K\setminus\{0\}$. The following conditions are equivalent:
\begin{enumerate}
\item There exist $\beta_1,\dots,\beta_d\in K$ such that
$$
\beta_1\alpha_1^n+\cdots+\beta_n\alpha_d^n=0
$$
for infinitely many $n\in \N$;
\item There exist positive integers $u,v\in \N$ and
$\gamma_1,\dots,\gamma_d\in \overline{K}$ such that
$\alpha_i=\gamma_i^u$ for all $i$ and
$\gamma_1^v,\dots,\gamma_d^v$ are linearly dependent over $\overline{\F}_p$;
\item If $L=K\cap \overline{\F}_p$, then there exist positive integers $u,v$
and elements $\mu_1,\dots,\mu_d\in L$ and $\gamma_1,\dots,\gamma_d\in
\overline{K}$, such that $\alpha_i=\mu_i\gamma_i^u$ for all $i$
and $\gamma_1^v,\dots,\gamma_d^v$ are linearly dependent over $L$.
\end{enumerate}
\end{theorem}
For example, the implication $(1)\Rightarrow (2)$ follows from the results
in this paper as follows. Suppose that
$$
a(n)=\beta_1\alpha_1^n+\cdots+\beta_n\alpha_d^n
$$
satisfies $a(n)=0$ for infinitely many $n$. If $a$ is not nondegenerate,
then $\alpha_i/\alpha_j$ is a root of unity for some $i\neq j$, 
say $(\alpha_i/\alpha_j)^N=1$.
Then clearly we may take $u=1$ and $v=N$ and $\gamma_i=\alpha_i$
for all $i$.
Thus, suppose that $a(n)$ is nondegenerate. By Theorem~\ref{thm1c}
 there exist $r,s\in \N$ and $p$-power $q$ 
 such that $a(r+sq^n)=0$ for all $n$. If we follow
 the proof of Lemma~\ref{lem5.8} then (2) follows. (In particular
 $x=0$ implies that $\alpha_1,\dots,\alpha_d$ are linearly dependent over
 $\F_q$).
 
It may be possible to use Masser's theorem to prove
 some of the results in this paper.

\section{Recurrence sequences in modules}
In this section we will study recurrence sequences in modules
over arbitrary commutative algebras. Our main goal is to
prove Theorem~\ref{thm1x} and Theorem~\ref{thm1b}. First we  prove that the intersection of two $p$-normal
sets is again $p$-normal.
\begin{definition}
A subset $N\subseteq \Z^r$ is called a {\it rectangular coset in $\Z^r$} if
it is of the form
$$
N=e_0+e_1\Z+e_2\Z+\cdots+e_s\Z
$$
where $e_0,e_1,e_2,\dots,e_s\in \Z^r$ and $e_1,\dots,e_s$
are pairwise orthogonal nonzero idempotent elements in $\Z^r$.
\end{definition}
Note that idempotent vectors in $\Z^r$ consist of $0$'s and $1$'s.
If $e_i$ and $e_j$ are orthogonal (i.e., $e_ie_j=0$), then there
is no position where 
$e_i$ and $e_j$ both have a $1$.
\begin{definition}
We call $N\subseteq \N^r$ a {\it rectangular coset in $\N^r$}
if it is of the form
$$
N=e_0+e_1\N+\cdots+e_s\N
$$
where $e_0,\dots,e_s\in \N^r$ and
$e_1,\dots,e_s$ are pairwise orthogonal nonzero idempotents.
\end{definition}
\begin{lemma}\label{lem9.2}
Suppose that $c_1,\dots,c_r\in \Q\setminus\{0\}$ and $q\in \Q$ with $q>1$.
The set
$$
A=\{(l_1,\dots,l_r)\in \Z^r\mid c_1q^{l_1}+\cdots+c_rq^{l_r}=0\}
$$
is a finite union of rectangular cosets in $\Z^r$.
\end{lemma}
\begin{proof}
We prove this lemma by induction on $r$. The case $r=1$ is clear.
 Choose $D\in \N$ such that
\begin{equation}\label{eqD}
|c_iq^{D}|>\sum_{j\neq i}|c_j|
\end{equation}
for all $i$. 

We claim that for every $(l_1,\dots,l_r)\in A$
there exist distinct indices $i$ and $j$ such that $l_j\leq l_i\leq D+l_j$: 
Take $i$ such that $l_i$ is maximal.
Suppose that $l_i>D+l_j$ for all $j\neq i$.
 Then we get
$$
0=\sum_{j=1}^r c_jq^{l_j},
$$
so
$$
|c_i|q^{l_i}=|c_iq^{l_i}|=|\sum_{j\neq i}c_jq^{l^j}|\leq \sum_{j\neq i}
|c_jq^{l_j}|\leq \sum_{j\neq i} |c_j|q^{l_i-D}
$$
and
$$
|c_i|q^{D}>\sum_{j\neq i}|c_j|.
$$
This is in contradiction with (\ref{eqD}).
Therefore, $l_i\leq D+l_j$ for some $j\neq i$.

We can write
$$
A=\bigcup_{i\neq j} \bigcup_{k=0}^D A_{i,j,k}
$$
where
$$
A_{i,j,k}=
A\cap \{(l_1,\dots,l_r)\mid l_i=l_j+k\}.
$$
Define
\begin{multline*}
B_{i,j,k}
=\big\{(l_1,\dots,l_{i-1},l_{i+1},\dots ,l_r)\in \Z^r\big| c_1q^{l_1}+\cdots+
c_{i-1}q^{l_{i-1}}+c_{i+1}q^{l_{i+1}}+\cdots\\
\cdots+(c_iq^k+c_j)q^{l_j}+\cdots+
c_rq^{l_r}=0\big\},
\end{multline*}
so that
$$
A_{i,j,k}=\{(l_1,\dots,l_{i-1},l_j+k,l_{i+1},\dots,l_r)\in \Z^{r-1}\mid 
(l_1,\dots,l_{i-1},l_{i+1},\dots,l_r)\in B_{i,j,k}\}.
$$
By the induction hypothesis, $B_{i,j,k}$ is a finite union of
rectangular cosets. Therefore, $A_{i,j,k}$ is a finite union of
rectangular cosets. We conclude that $A$ is a finite union of rectangular
cosets.
\end{proof}
We also have the following monoid version of this Lemma~\ref{lem9.2}.
\begin{corollary}\label{correct}
Suppose that $c_1,\dots,c_r\in \Q$ are nonzero and $q\in \Q$ with $q>1$.
The set
$$
A=\{(l_1,\dots,l_r)\in \N^r\mid c_1q^{l_1}+\cdots+c_rq^{l_r}=0\}
$$
is a finite union of rectangular cosets of in $\N^r$.
\end{corollary}
\begin{proof}
This follows from the previous lemma and the observation
that a rectangular coset in $\Z^r$ intersected with $\N^r$ is
a rectangular coset in $\N^r$.
\end{proof}

\begin{lemma}\label{lem8.2}
Suppose that $S_1,S_2\subseteq\N$.
are both $p$-normal of order $\leq d$, then
 $S_1\cap S_2$ is $p$-normal of order $\leq d$.
\end{lemma}
\begin{proof}
Without loss of generality we may assume that $S_1$ is
an infinite arithmetic progression or an elementary $p$-nested set
of order $\leq d$.
Similarly, we may assume that $S_2$ is an
infinite arithmetic progression or an elementary $p$-nested set
of order $\leq d$.

{\bf case 1:} The lemma is clear if 
both $S_1$ and $S_2$ are infinite arithmetic progressions.

{\bf case 2:} Suppose $S_1=m+n\N$  is an infinite arithmetic progression
and $S_2=S_q(c_0;c_1,\dots,c_r)$ with $r\leq d$. Without loss of generality
we may assume that $m<n$. Then we have $S_1=(m+n\Z)\cap \N$. 
We can write $n=p^lu$
with $l\in \N$ and $\gcd(u,p)=\gcd(u,q)=1$. 
Since $m+n\N=(m+p^l\N)\cap (m+u\N)$ by 
the Chinese Remainder Theorem, we may assume that either $\gcd(n,q)=1$
or $n$ divides some power of $q$.

{\bf case 2a:} Suppose that $\gcd(n,q)=1$.
Choose $s$ such that $q^s\equiv 1 \bmod n(q-1)$.
There exists a decomposition
$$
S_q(c_0;c_1,\dots,c_r)=\bigcup_{0\leq l_1,\dots,l_r<s} 
S_{q^s}(c_0;c_1q^{l_1},\dots,c_rq^{l_r}).
$$
There is an inclusion
$$
S_{q^s}(c_0;c_1q^{l_1},\dots,c_rq^{l_r})\subseteq 
c_0+c_1q^{l_1}+\cdots c_rq^{l_r}+n\Z,
$$
because $(q-1)c_i\in \Z$ for all $i$ and $q^s\equiv 1 \bmod (q-1)n$.
It follows that $S_1\cap S_2=S_q(c_0;c_1,\dots,c_r)\cap (m+n\Z)$ is the union
of all $S_{q^s}(c_0;c_1q^{l_1},\dots,c_rq^{l_r})$ for which
$$
c_0+c_1q^{l_1}+\cdots+c_rq^{l_r}\equiv m\bmod n.
$$
So $S_1\cap S_2$ is $p$-normal of order $\leq d$.

{\bf Case 2b:} Suppose that $n$ divides $q^s$ for some positive integer $s$. Then $S_q(c_0;c_1,\dots,c_r)$
is the union  of all
$$
S_q(c_0+c_{i_1}q^{l_1}+\cdots,c_{i_u}q^{l_u},c_{j_1}q^s,c_{j_2}q^s,
\dots,c_{j_v}q^s)
$$
for which $\{1,2,\dots,r\}$ is a disjoint union of $\{i_1,\dots,i_u\}$
and $\{j_1,\dots,j_v\}$, $r=u+v$ and $0\leq l_1,l_2,\dots,l_u<s$ .
Note that
$$
S_q(c_0+c_{u_1}q^{l_1}+\cdots,c_{i_u}q^{l_u},c_{j_1}q^s,c_{j_2}q^s,
\dots,c_{j_v}q^s)\subseteq 
c_0+c_{u_1}q^{l_1}+\cdots,c_{i_u}q^{l_u}+n\Z.
$$
because $n$ divides $q$. We conclude that $S_1\cap S_2$ is
$p$-normal of order $\leq d$ as in case 2a.

{\bf case 3:}
Suppose that $S_1=S_q(c_0;c_1,\dots,c_r)$ and
$S_2=S_{q'}(c_0';c_1',\dots,c_{r'}')$.
If $q''$ is an integral power of $q$ and also an integral power of $q'$
then both $S_1$ and $S_2$ can be written as a finite union
of sets of the form $S_{q''}(f_0;f_1,\dots,f_s)$.
We can reduce to the case where $q'=q$.

The set 
$$
M=\big\{(l_1,\dots,l_r,l_1',\dots,l_{r'}')\in \N^{r+r'}\big| 
c_0+c_1q^{l_1}+\cdots+c_rq^{l_r}=
e_0+e_1q^{l_1'}+\cdots+e_{r'}q^{l_{r'}'}\big\}
$$
is a finite union of rectangular cosets (see Corollary~\ref{correct}). 
From this it follows that
$$
\big\{(l_1,\dots,l_r)\in \N^r\big|\exists l_1',\dots,l_{r'}'\in \N\,
c_0+c_1q^{l_1}+\cdots+c_rq^{l_r}=
e_0+e_1q^{l_1'}+\cdots+e_{r'}q^{l_{r'}'}\big\}
$$
is also a finite union of rectangular cosets in $\N^r$. 
Therefore
$$
S_q(c_0;c_1,\dots,c_r)\cap S_q(c_0';c_1',\dots,c_{r'}')
$$
is a finite union of sets of the form
$$
S_{q}(f_0;f_1,\dots,f_u)
$$
with $u\leq r\leq d$.
This shows that
$$
S_1\cap S_2=S_q(c_0;c_1,\dots,c_r)\cap S_q(c_0';c_1',\dots,c_{r'}')
$$
is $p$-normal of order $\leq d$.
\end{proof}

Suppose $R$ is a commutative ring, $M$ is an $R$-module and 
${\mathfrak p}\subset R$ is a
prime ideal. A submodule $N\subseteq M$ is called ${\mathfrak p}$-primary
if ${\mathfrak p}$ is the only associated prime of $M/N$. We call
$M$ ${\mathfrak p}$-coprimary if the only associated prime of $M$ is
${\mathfrak p}$ (i.e., if the submodule $(0)\subseteq M$ is ${\mathfrak
p}$-primary). (See \cite[\S3.3]{Eis})
\begin{lemma}\label{lemprimary}
Suppose that $K$ is a field, $R$ is a finitely generated $K$-algebra and
 $M$ is a finitely generated $R$-module.
If $a\in M^{\N}$ is an $R$-recurrence sequence of order $d$, then
there exist prime ideals ${\mathfrak p}_1,\dots,{\mathfrak p}_l\in R$,
finitely generated $R$-modules $M_1,\dots,M_l$ where $M_i$ is ${\mathfrak p}_i$-coprimary
 for all $i$ and  $R$-recurrence sequences $a_1,a_2,\dots,a_l$
of order $\leq d$
with $a_i\in M_i^{\N}$ such that
$$
{\mathcal Z}(a)=\bigcap_{i=1}^l{\mathcal Z}(a_i).
$$
\end{lemma}
\begin{proof}
We use the primary decomposition of $M$.
There exist prime ideals ${\mathfrak p}_1,\dots,{\mathfrak p}_l$
and submodules $N_1,N_2,\dots,N_l$ of $M$ such that
$\bigcap_{i=1}^lN_i=0$ and $N_i$ is a ${\mathfrak p}_i$-primary submodule
of $M$
for $i=1,2,\dots,l$ (see~\cite[Theorem 3.10]{Eis}).
Let $\rho_i:M\to M_i:=M/N_i$ be the quotient homomorphism.
Because the intersection of the kernels of $\rho_i,i=1,2,\dots,l$
is equal to $0$, we get
$$
{\mathcal Z}(a)=\bigcap_{i=1}^l{\mathcal Z}(\rho_i(a)).
$$
Now take $a_i=\rho_i(a)\in M_i^{\N}$.
\end{proof}
\begin{lemma}\label{lemreducefield}
Suppose that $K$ is a  field of characteristic $p$, 
$R$ is a finitely domain over $K$
and $M$ is
a finitely generated torsion-free $R$-module. Suppose that $a\in M^{\N}$
is an $R$-recurrence sequence of order $d$.
If $p=0$ then ${\mathcal Z}(a)$ is a union of a finite set
and finitely many infinite arithmetic progressions.
If $p>0$,
 then
${\mathcal Z}(a)$ is $p$-normal of order $\leq d-2$.
\end{lemma}
\begin{proof}
Let $L$ be the quotient field of $R$. We may view $a$ as
an $L$-recurrence sequence of order $\leq d$ in the vector space 
$V=M\otimes_RL\supseteq M$. Choose a basis $e_1,\dots,e_r$ of $V$.
We can write $a=\sum_{i=1}^ra_ie_i$ with $a_i\in L^{\N}$ a
$L$-recurrence sequence of order $\leq d$ for all $i$. 
We have
$$
{\mathcal Z}(a)=\bigcap_{i=1}^r {\mathcal Z}(a_i).
$$
If $p=0$ then each ${\mathcal Z}(a_i)$ is a union
of a finite set and finitely many arithmetic progressions
by the Skolem-Mahler-Lech theorem (Theorem~\ref{theoSML}).
Hence, ${\mathcal Z}(a)$ is a union of a finite set
and finitely many infinite arithmetic progressions.

If $p>0$, then ${\mathcal Z}(a_i)$ is $p$-normal of order $\leq d-2$ for all 
$i$ by Theorem~\ref{thm1a}.
This implies that ${\mathcal Z}(a)$ is $p$-normal of order $\leq d-2$
by Lemma~\ref{lem8.2}.
\end{proof}

\begin{lemma}\label{lem10}
Suppose that $K$ is an infinite field of characteristic $p$,
$R$ is a finitely generated
$K$-algebra, ${\mathfrak p}\subseteq R$ is a prime ideal,
$M$ is a finitely generated ${\mathfrak p}$-coprimary module,
and $a\in M^{\N}$ is a linear $R$-recurrence sequence.
If $p=0$ then ${\mathcal Z}(a)$ is a union of a finite set
and finitely arithmetic progressions.
If $p>0$ then
${\mathcal Z}(a)$ is $p$-normal.
\end{lemma}
\begin{proof}
By the Noether Normalization Lemma (see for example~\cite[\S8.2.1]{Eis}), 
there exist algebraically independent
$x_1,x_2,\dots,x_s\in R/{\mathfrak p}$ such that
$R/{\mathfrak p}$ is a finite 
$K[x_1,\dots,x_s]$-module.
Choose $y_1,\dots,y_s\in R$
such that $y_i+{\mathfrak p}=x_i$ for all $i$.

Now $M$ is also a $K[y_1,\dots,y_s]$-module.
Since $M$ is ${\mathfrak p}$-primary, there exists
$t\in \N$ such that ${\mathfrak p}^tM=0$ (see~\cite[Proposition~3.9]{Eis}).
We have a filtration
$$
M\supseteq {\mathfrak p}M\supseteq {\mathfrak p}^2M\supseteq
\cdots \supseteq {\mathfrak p}^tM=0.
$$
and each quotient ${\mathfrak p}^iM/{\mathfrak p}^{i+1}M$
is a finite $R/{\mathfrak p}R$-module, hence a finite
$K[y_1,\dots,y_s]$-module. 
It follows that $M$  is a finite
$K[y_1,\dots,y_s]$-module. 
Similarly, $R/{\mathfrak p}^tR$ is a finite
$K[y_1,\dots,y_s]$-module.
The $R$-module generated by $a,Ea,E^2a,\dots$ is finitely generated
as an $R/{\mathfrak p}^t$-module, hence
it is  finitely generated as a
$K[y_1,\dots,y_s]$-module. Therefore, $a$ satisfies a
$K[y_1,\dots,y_s]$-recurrence relation. However, the
order of $a$ as a $K[y_1,\dots,y_s]$-recurrence sequence can be larger
than the order of $a$ viewed as an $R$-recurrence sequence.
Since ${\mathfrak p}\cap K[y_1,\dots,y_s]=(0)$, 
the annihilator in $K[y_1,\dots,y_s]$ of any nonzero element in $M$
is the zero ideal. 
We can reduce to the situation where
$R=K[y_1,\dots,y_s]$ is the polynomial ring and $M$ is a finitely generated
$R$-module without torsion. We now apply Lemma~\ref{lemreducefield}
\end{proof}
\begin{lemma}\label{lemfg}
Suppose that $K$ is a field, $R$ is a $K$-algebra 
and $a\in M^{\N}$ is a linear $R$-recurrence
sequence of order $d$. Then there exists a ring $R'$ which
is finitely generated over $K$,
a finitely generated $R'$-module $M'$ and a $R'$-recurrence
sequence $a'\in (M')^{\N}$ such that
$$
{\mathcal Z}(a)={\mathcal Z}(a').
$$
\end{lemma}
\begin{proof}
There exists a nonnegative integer $m$ and
$\alpha_0,\dots,\alpha_{m-1}\in R$ such that
$$
E^ma=\sum_{i=0}^{m-1}\alpha_iE^ia.
$$
Set $R':=K[\alpha_0,\dots,\alpha_{m-1}]$.
Then $a$ is also an linear $R'$-recurrence sequence.
Let $N$ be the $R'$-module generated by $a,Ea,E^2a,E^3a,\dots$.
Then $N$ is finitely generated. In fact, it is generated
by $a,Ea,E^2a,\dots,E^{m-1}a$.
Let $M'$ be the module generated by $a(0),a(1),a(2),\dots$.
Then $M'$ is a finitely generated $R'$-module because it
is a homomorphic image of the finitely generated module $N$
via the homomorphism $a\mapsto a(0)$. Take $a'=a$.
\end{proof}
\begin{proof}[Proof of Theorem~\ref{thm1x}]
 By Lemma~\ref{lemfg}
 we may assume that $R$ is finitely generated $\Q$-algebra and $M$
 is finitely generated as an $R$-module.
We apply  Lemma~\ref{lemprimary}. There exist prime ideals ${\mathfrak p}_1,
\dots,{\mathfrak p}_l\subseteq R$, $R$-modules
$M_1,M_2,\dots,M_l$ with $M_i$ ${\mathfrak p}_i$-coprimary for all $i$,
recurrence sequences $a_i\in M_i^{\N}$ for all $i$  such that
$$
{\mathcal Z}(a)=\bigcap_{i=1}^l {\mathcal Z}(a_i).
$$
By Lemma~\ref{lem10}, ${\mathcal Z}(a_i)$ is a union of a finite set 
and finitely many arithmetic progressions for all $i$. But then
${\mathcal Z}(a)$ also is a union of a finite set and finitely
many arithmetic progressions.
\end{proof}

\begin{lemma}\label{lemdensity}
Suppose that $A$ is an elementary $p$-nested set
of order $l$   which is up to a finite set 
contained in a $p$-nested set $B$ of order $d$.
Then $l\leq d$.
\end{lemma}
\begin{proof}
This follows from 
$$
\delta_A(n)=\Omega(\log(n)^l)
$$
and
$$
\delta_B(n)=O(\log(n)^d).
$$
\end{proof}

\begin{proof}[Proof of Theorem~\ref{thm1b}]
Suppose that $K$ is an infinite field containing $\F_p$ (for example
$K=\overline{\F}_p$, the algebraic closure). We may view $a$ as
a $R\otimes_{{\F}_p}K$-sequence in $M\otimes_{\F_p}K\supseteq M$ rather than in $M$.
This shows that we may assume that $R$ contains
an infinite field without loss of generality.

By Lemma~\ref{lemfg}
we may assume that $R$ is finitely generated $K$-algebra and $M$
 is finitely generated as an $R$-module.
By Lemma~\ref{lemprimary} and Lemma~\ref{lem8.2},
we can reduce to the case where there exists
a prime ideal ${\mathfrak p}\subseteq R$ such that $M$ is ${\mathfrak p}$-coprimary.

By Lemma~\ref{lem10}, ${\mathcal Z}(a)$ is $p$-normal. However,
it is not yet clear that ${\mathcal Z}(a)$ is $p$-normal {\it of order 
$\leq d-2$}. We prove Theorem~\ref{thm1b} by induction on $t$
where $t$ is the smallest nonnegative integer such that ${\mathfrak p}^tM=0$.
The case $t=0$ is clear.

We can write
$$
{\mathcal Z}(a)=\big(F_1\cup \bigcup_{i=1}^{\alpha} A_i\cup 
\bigcup_{i=1}^{\beta} B_i\big)\setminus F_2
$$
where $F_1,F_2$ are finite, $A_i$ is an infinite arithmetic
progression for all $i$ and $B_i$ elementary $p$-nested
for all $i$. We may assume that all the arithmetic progressions $A_i$ have the
same period, say $n$. Since the intersection
of a set  of the form $S_q(c_0;c_1,\dots,c_l)$ with $m+n\Z$ is
a union of elementary $p$-nested sets of order $\leq l$
(see the proof of Lemma~\ref{lem8.2}, case 2), we
may assume that each $B_j$ is contained in $m+n\Z$ for some $m\in \Z$.
Hence we may assume that $B_j\cap A_i=\emptyset$
for all $i,j$.

Suppose that $B_j=S_q(c_0;c_1,\dots,c_l)$. We would like
to show that $l\leq d-2$. 
Let us write ${\mathfrak p}=(f_1,\dots,f_k)$. 
Define 
$$a'\in ({\mathfrak p}M)^k)^{\N}\cong (({\mathfrak p}M)^{\N})^k
$$
by
$$
a'=(f_1a,f_2a,\dots,f_ka).
$$
Then $a'$ is a $R$-recurrence sequence of order $\leq d$.
By induction on $t$, ${\mathcal Z}(a')$ is $p$-normal of order $\leq d-2$.
We  can write
$$
{\mathcal Z}(a')=(\bigcup_{i=1}^{\alpha'} A_i'\cup B')\setminus F',
$$
where $F'$ is finite, $A_i'$ is an infinite arithmetic
progression for all $i$ and $B$ is a $p$-nested set of
order $d-2$. 

If $B_j\cap A_i'$ is finite for all $i$, then $B_j$ is up to
a finite set contained
in $B'$. We get $l\leq d-2$ by Lemma~\ref{lemdensity}.

Suppose that $B_j\cap A_i'$ is
infinite for some $i$. We can write $A_i'=m'+n'\N$ for some $m',n'\in \Z$.
Then there exists a $p$-power $q'$ and $c_0',\dots,c_l'\in \Q$ such
that 
$$S_{q'}(c_0';c_1',\dots,c_l')\subseteq B_j\cap A_i'\subset {\mathcal Z}(a)\cap
(m'+n'\Z).$$
Now we have
$$
{\mathcal Z}(a)\cap (m'+n'\Z)=L_{m'}^{n'}({\mathcal Z}(T_{m'}^{n'}a)).
$$
Since $m'+n'i\in {\mathcal Z}(a')$ for $i\gg 0$, we get
$$
(T_{m'}^{n'}a)(i)\in 
N:=\{f\in N\mid {\mathfrak p}f=(0)\}
$$
for $i\gg 0$. Note that $N$ is a finitely generated 
torsion-free $R/{\mathfrak p}$-module, 
and $T_{m'}^{n'}a$ is an $R/{\mathfrak p}$-recurrence
sequence of order $\leq d$.
By Lemma~\ref{lemreducefield}, ${\mathcal Z}(T_{m'}^{n'}a)$ is $p$-normal of 
order $\leq d-2$.
We have
$$
S_q(c_0';c_1',\dots,c_l')\subseteq L^{n'}_{m'}({\mathcal
Z}(T^{n'}_{m'}a))\subseteq {\mathcal Z}(a)$$
and $L^{n'}_{m'}({\mathcal Z}(T^{n'}_{m'}a))$ is $p$-normal of 
order $\leq d-2$. Since $S_q(c_0';c_1',\dots,c_l')\subseteq B_j$
and $B_j\cap A_k=\emptyset$ for all $k$, the intersection
of $S_q(c_0';c_1',\dots,c_l')$
and any infinite arithmetic progression in ${\mathcal Z}(a)$ is finite.
In particular, the intersection of $S_q(c_0';c_1',\dots,c_l')$
and any infinite arithmetic progression
in $L^{n'}_{m'}({\mathcal Z}(T^{n'}_{m'}a))$ is finite.
It follows that $S_q(c_0';c_1'\dots,c_l')$ is up to a finite
set contained in a $p$-nested set of order $\leq d-2$, so
$l\leq d-2$ by Lemma~\ref{lemdensity}.
\end{proof}

\section{Open problems}
{\bf Uniform bounds.} Suppose that $K$ is a field of characteristic
$p>0$ and $a\in K^\N$ is a recurrence sequence of order $d$ such 
that ${\mathcal Z}(a)$ is finite. Does there exist
a uniform upper bound for $|{\mathcal Z}(a)|$ depending
only on $d$ and $K$? In characteristic 0 such uniform
bound exists for number fields (see for example~\cite{ES1,ES2,ESS,ScSc}).

{\bf $S$-unit equations.}
Suppose that $K$ is a field and $\Gamma\subseteq K^\star$ is a 
finitely generated subgroup. Fix $\beta_1,\dots,\beta_d\in K$
and consider the set ${\mathcal F}_{\Gamma}(\beta_1,\dots,\beta_d)$
of all $(\alpha_1,\dots,\alpha_d)\in \Gamma^d$ such that
\begin{equation}\label{eq111}
\beta_1\alpha_1+\cdots+\beta_d\alpha_d=1
\end{equation}
and no proper subsum 
$$
\beta_{i_1}\alpha_{i_1}+\cdots+\beta_{i_m}\alpha_{i_m}
$$
with $1\leq i_1<i_2<\cdots<i_m\leq d$ and $1\leq m<d$ is equal to $0$.
If $K$ has characteristic 0 then  
${\mathcal F}_{\Gamma}(\beta_1,\dots,\beta_d)$
is finite (see~\cite{VDP}). 
In that case an upper bound for the number of elements
of ${\mathcal F}_{\Gamma}(\beta_1,\dots,\beta_d)$ may
be given in terms of the number of generators of $\Gamma$ and $d$
(see~\cite{ESS,ES1,ES2}).

If $K$ has positive characteristic then ${\mathcal
F}_{\Gamma}(\beta_1,\dots,\beta_d)$ may be infinite. Solutions of (\ref{eq111})
in positive characteristic were studied by Masser in \cite{Masser}.
His results are sufficient to solve
a conjecture of Klaus~Schmidt about mixing properties of algebraic $\Z^{r}$-actions.
Our methods here may lead to  a precise
description of the solutions of (\ref{eq111}) in positive characteristic.

{\bf A nonlinear Skolem-Mahler-Lech theorem.}
A set of the form $m+n\Z\subseteq \Z$ with $m\in \Z$ and $n$
a positive integer is called a doubly infinite arithmetic progression.
The following theorem was recently proven by Jason Bell (\cite{Bell}).
\begin{theorem}
Suppose that $Y$ is an affine variety over a field $K$ of
characteristic 0, $X\subseteq Y$ is a Zariski closed subset, $y\in Y$, and
 $\sigma:Y\to Y$ is an automorphism.
 The set
 $$
 \{n\in\Z\mid \sigma^n(y)\in X\}
 $$
 is a union of a finite set and a finite number
 of doubly infinite arithmetic progressions.
 \end{theorem}
 Bell generalized the methods used in the proof
 of the Skolem-Mahler-Lech theorem (in characteristic 0) to prove this result.
Given the results in this paper, it is natural to
conjecture the following:
\begin{conjecture}
Suppose that $Y$ is an affine variety over a field $K$ of characteristic
$p>0$, $X$ is a Zariski closed subset, $y\in Y$ and
$\sigma:Y\to Y$ is an automorphism.
Then the set
$$
\{n\in \Z\mid \sigma^n(y)\in X\}
$$
is a union of a finite set, finitely many doubly infinite arithmetic
progressions and finitely many sets of the form
$$
\widetilde{S}_q(c_0;c_1,c_2,\dots,c_s).
$$
\end{conjecture}

{\sl
Harm Derksen\\
Department of Mathematics\\
University of Michigan\\
hderksen@umich.edu}

\end{document}